\documentclass[12pt,reqno]{amsart}
\usepackage[usenames]{color}
\usepackage[francais]{babel} 
 \usepackage[applemac]{inputenc}
\usepackage{graphicx}
\usepackage{amsmath,epsf}
\usepackage[colorlinks=true,citecolor=cyan]{hyperref}
\usepackage{epsfig}
\usepackage{amsfonts}
\usepackage{amssymb}
\usepackage{bbm}
\usepackage{amsthm}
\usepackage[mathscr]{eucal}
\usepackage{wrapfig}

\def\longuefleche#1#2{\put(#1,#2){\rouge{\put(0,.19){\line(1,0){3}}\put(2.3,0){$\rightarrow$}}}}
\def\flecheverslebas#1#2{\put(#1,#2){\rouge{\put(0,0){\line(0,-1){2}}\put(-0.19,-1.8){$\downarrow$}}}}
\def\flecheverslehaut#1#2{\put(#1,#2){\rouge{\put(0,0){\line(0,-1){2}}\put(-0.19,-0.5){$\uparrow$}}}}
\def\pointille#1#2{\put(#1,#2){\bleu{\put(0,0){\bf\ --\ --\ --\ }\put(2.5,0.02){$\rightarrow$}}}}

\def\defn#1{{\bleu{\bf #1}}}

\def\ligne(#1,#2)#3{{\multiput(#1,#2)(1,0){#3}{\carre}}}
\def\jligne(#1,#2)#3{{\multiput(#1,#2)(1,0){#3}{\jcarre}}}
\def\colonne(#1,#2)#3{{\multiput(#1,#2)(0,1){#3}{\carre}}}
\def\jcolonne(#1,#2)#3{{\multiput(#1,#2)(0,1){#3}{\jcarre}}}

\def\seriequatre#1#2#3#4#5#6{\bleu{#1}+\bleu{#2}\,#6+\bleu{#3}\,\frac{#6^2}{2}+\bleu{#4}\,\frac{#6^3}{6}+\bleu{#5}\,\frac{#6^4}{24}}
\def\seriecinq#1#2#3#4#5#6#7{\bleu{#1}+\bleu{#2}\,#7+\bleu{#3}\,\frac{#7^2}{2}+\bleu{#4}\,\frac{#7^3}{6}+\bleu{#5}\,\frac{#7^4}{24}+\bleu{#6}\,\frac{#7^5}{120}+\ \ldots}

\def\card#1{\left|#1\right|}

\newdimen\graphelength
\def\grapheun#1#2#3{\setlength{\unitlength}{\graphelength}
\begin{picture}(73,40)(-5,22)
\put(-4,-5){\framebox(70,65)}
\put(15,7){\bleu{\circle*{7}}}  
        \put(0,0){$#2$}   
\put(45,7){\bleu{\circle*{7}}}
        \put(52,0){$#3$}
\put(30,37){\bleu{\circle*{7}}}
        \put(25,43){$#1$}
\end{picture}}

\def\graphedeux#1#2#3{\setlength{\unitlength}{\graphelength}
\begin{picture}(73,40)(-5,22)
\put(-4,-5){\framebox(70,65)}
\put(15,8){\rouge{\line(1,2){15}}}
\put(15,7){\bleu{\circle*{7}}}  
        \put(0,0){$#2$}   
\put(45,7){\bleu{\circle*{7}}}
        \put(52,0){$#3$}
\put(30,37){\bleu{\circle*{7}}}
        \put(25,43){$#1$}
\end{picture}}

\def\graphetrois#1#2#3{\setlength{\unitlength}{\graphelength}
\begin{picture}(73,40)(-5,22)
\put(-4,-5){\framebox(70,65)}
\put(45,8){\rouge{\line(-1,2){15}}}
\put(15,7){\bleu{\circle*{7}}}  
        \put(0,0){$#2$}   
\put(45,7){\bleu{\circle*{7}}}
        \put(52,0){$#3$}
\put(30,37){\bleu{\circle*{7}}}
        \put(25,43){$#1$}
\end{picture}}

\def\graphequatre#1#2#3{\setlength{\unitlength}{\graphelength}
\begin{picture}(73,40)(-5,22)
\put(-4,-5){\framebox(70,65)}
\put(45,7){\rouge{\line(-1,0){30}}}
\put(15,7){\bleu{\circle*{7}}}  
        \put(0,0){$#2$}   
\put(45,7){\bleu{\circle*{7}}}
        \put(52,0){$#3$}
\put(30,37){\bleu{\circle*{7}}}
        \put(25,43){$#1$}
\end{picture}}

\def\graphecinq#1#2#3{\setlength{\unitlength}{\graphelength}
\begin{picture}(73,40)(-5,22)
\put(-4,-5){\framebox(70,65)}
\put(15,8){\rouge{\line(1,2){15}}}
\put(45,8){\rouge{\line(-1,2){15}}}
\put(15,7){\bleu{\circle*{7}}}  
        \put(0,0){$#2$}   
\put(45,7){\bleu{\circle*{7}}}
        \put(52,0){$#3$}
\put(30,37){\bleu{\circle*{7}}}
        \put(25,43){$#1$}
\end{picture}}

\def\graphesix#1#2#3{\setlength{\unitlength}{\graphelength}
\begin{picture}(73,40)(-5,22)
\put(-4,-5){\framebox(70,65)}
\put(15,8){\rouge{\line(1,2){15}}}
\put(45,7){\rouge{\line(-1,0){30}}}
\put(15,7){\bleu{\circle*{7}}}  
        \put(0,0){$#2$}   
\put(45,7){\bleu{\circle*{7}}}
        \put(52,0){$#3$}
\put(30,37){\bleu{\circle*{7}}}
        \put(25,43){$#1$}
\end{picture}}

\def\graphesept#1#2#3{\setlength{\unitlength}{\graphelength}
\begin{picture}(73,40)(-5,22)
\put(-4,-5){\framebox(70,65)}
\put(45,8){\rouge{\line(-1,2){15}}}
\put(45,7){\rouge{\line(-1,0){30}}}
\put(15,7){\bleu{\circle*{7}}}  
        \put(0,0){$#2$}   
\put(45,7){\bleu{\circle*{7}}}
        \put(52,0){$#3$}
\put(30,37){\bleu{\circle*{7}}}
        \put(25,43){$#1$}
\end{picture}}

\def\graphehuit#1#2#3{\setlength{\unitlength}{\graphelength}
\begin{picture}(73,40)(-5,22)
\put(-4,-5){\framebox(70,65)}
\put(15,8){\rouge{\line(1,2){15}}}
\put(45,8){\rouge{\line(-1,2){15}}}
\put(45,7){\rouge{\line(-1,0){30}}}
\put(15,7){\bleu{\circle*{7}}}  
        \put(0,0){$#2$}   
\put(45,7){\bleu{\circle*{7}}}
        \put(52,0){$#3$}
\put(30,37){\bleu{\circle*{7}}}
        \put(25,43){$#1$}
\end{picture}}
\newcommand{\pref}[1]{(\ref{#1})}

\def\Id{{\rm Id}}
\def\bleu{\textcolor{blue}}
\def\rouge{\textcolor{red}}

\def\blanc{\textcolor{white}}

\def\ssi{\quad {\rm ssi}\quad}

\def\P{\bleu{\mathcal{P}}}

\def\Gra{\bleu{\mathsf{Gra}}}
\def\Gro{\bleu{\mathsf{Gro}}}

\def\Part{\bleu{\mathsf{Part}}}

\def\Arb{\bleu{\mathsf{Arb}}}
\def\Arbo{\bleu{\mathscr{A}}}
\def\Vertebre{\bleu{\mathscr{V}}}
\def\Endo{\bleu{\mathsf{End}}}

\def\Cyc{\bleu{\mathcal C}}
\def\Der{\bleu{\mathsf{Der}}}

\def\Bin{\bleu{\mathscr{B}}}
\def\Ens{\bleu{\mathsf{E}}}
\def\Inv{\bleu{\mathsf{Inv}}}
\def\SingletonX{\bleu{\mathsf{X}}}

\def\Zero{\bleu{\mathsf{O}}}
\def\Un{\bleu{\mathbbm{1}}}
\def\Liste{\bleu{\mathbb{L}}}

\def\point{{^{\raise 0.7mm \hbox{$\scriptscriptstyle \bullet$}}}} 
\def\nonvide{^{\scriptscriptstyle +}}

\def\especeF{\bleu{F}}
\def\especeG{\bleu{G}}

\def\S{\bleu{{\mathbb S}}}

\def\fleche{\rightarrow}
\def\image#1{\,\rouge{\buildrel{#1}\over\mapsto}\,}
\def\bijection{\smash{\buildrel\sim\over\longrightarrow}}

\newtheorem{thm}{Théorème}[section]
\newtheorem{prop}[thm]{Proposition}

\theoremstyle{definition}

\title[Esp\`eces]{\bleu{Qu'est-ce qu'une esp\`ece de structures? \\Gen\`Ese et description}}

\author{F. Bergeron et G. Labelle}
\address{D\'epartement de Math\'ematiques, UQAM,  C.P. 8888, Succ. Centre-Ville, 
 Montr\'eal,  H3C 3P8, Canada.}

\begin{document}
\maketitle

 {
   \setcounter{tocdepth}{1}
   \parskip=0pt
   \footnotesize
   \tableofcontents
 }

\section{Introduction.}\label{sec_intro_especes}
Le concept de structure discrète est fondamental en science et se retrouve dans toutes les branches des mathématiques ainsi qu'en informatique théorique, où il porte le nom de structure de données.  L'un des objectifs de la théorie des espèces de structures est d'en donner une description claire. Cette définition se fait dans un esprit similaire à celui de la définition de fonction. Rappelons que la notion de fonction a subi une lente et obligatoire évolution, dans l'histoire des mathématiques, pour aboutir à la notion moderne. Cette vision moderne a séparé la définition de la notion de fonction du mode de description de fonctions explicites. Bien entendu, rien n'empêche ensuite de décrire une fonction particulière par une formule, un algorithme, comme solution d'une équation, etc. Là où l'on gagne, c'est dans la possibilité de définir de façon claire des espaces de fonctions. C'est un gain indéniable quand on pense à toutes les répercussions dans l'évolution de certains domaines des mathématiques (comme l'analyse fonctionnelle) ou de la physique quantique.

C'est de ce niveau de généralité que relève la définition des espèces de structures. Le but est donc de caractériser ce qu'est une structure combinatoire (structure discrète), plutôt que d'imposer une forme explicite de description. Accessoirement, la notion d'espèce permet aussi d'asseoir fermement un autre concept important: celui d'isomorphisme de structures combinatoires. Cette notion joue un rôle central dans la comparaison de structures, et permet de déduire une foule d'identités de façon très \og géométrique\fg. D'autres théorie ont été suggérées, mais aucune ne met aussi clairement  l'emphase sur ce
côté fonctionnel. 

\section{Genèse.}\label{sec_genese}
Au moins depuis Euler, l'utilisation de séries ou de polynômes formels pour démontrer des identités combinatoires s'est avérée particulièrement efficace, bien qu'un peu mystérieuse. Au départ, tout repose sur le fait qu'on a égalité entre séries
   $$\bleu{a_0}+\bleu{a_1}\,x+\bleu{a_2}\,x^2+\bleu{a_3}\,x^3+\ \ldots\ =
       \rouge{b_0}+\rouge{b_1}\,x+\rouge{b_2}\,x^2+\rouge{b_3}\,x^3+\ \ldots\,
          $$
  si et seulement si
     $$\bleu{a_0}= \rouge{b_0},\quad \bleu{a_1}= \rouge{b_1},\quad \bleu{a_2}= \rouge{b_2},\quad \bleu{a_3}= \rouge{b_3},\quad \ldots$$
  On exploite ensuite judicieusement (c'est là que réside l'aspect un peu mystérieux) les opérations entre séries pour obtenir des identités de la combinatoire énumérative.
Une illustration classique de cette efficacité est dans une dérivation simple de l'identité
\begin{equation}\label{id_binom}
   \bleu{ \binom{n}{k}=\binom{n-1}{k}+\binom{n-1}{k-1}}.
 \end{equation}
Rappelons que les \defn{coefficients binomiaux} apparaissent (par définition) comme coefficients du polynôme
 $$\bleu{ (x+1)^n=\sum_{k=0}^n \binom{n}{k}\, x^k}.$$
On constate alors que l'identité \pref{id_binom} découle  simplement de la comparaison du coefficient de $x^k$ dans chaque membre de l'égalité évidente: $\bleu{(x+1)^n=(x+1)\,(x+1)^{n-1}}$.
Dans des situations plus complexes, le fait de savoir trouver comment démontrer de cette façon des identités combinatoires est longtemps resté du grand art. Ce n'est qu'au début des années 1980 que des théories claires ont commencé à émerger pour rendre le tout limpide.

L'une des théorie les plus efficaces et englobantes, est la \defn{théorie des espèces de structures} proposée par André Joyal autour de 1980. 
Il a montré comment adapter la notion d'espèces de structures, de Charles Ehresmann (1905-1979), au contexte de la combinatoire énumérative. Plus particulièrement, en développant un lien clair entre opérations sur des espèces (de nature géométrique) et opérations entre séries formelles (plus algébrique), il a nettement mis en évidence comment dériver de façon limpide des identités difficiles à obtenir autrement. Très rapidement, une équipe de chercheurs de l'UQAM (F.B., P. Leroux, J. Labelle, et G.L., et leurs étudiants) s'est alors mise à explorer les ramifications de la théorie en en montrant la puissance et en en développant de nouveaux aspects.

\section{Espèces de structures.}\label{sec_especes}
La notion \defn{espèce de structures}\footnote{On peut trouver une description plus détaillée de la théorie des espèces et de ses applications dans la monographie \cite{BLL}.} $\especeF$  contient deux parties. 
\begin{enumerate}
\item[1)] 
Une première partie qui décrit comment  produire, pour chaque ensemble fini  $A$, un ensemble fini  $\especeF[A]$. On dit que les éléments $s$ de $\especeF[A]$ sont les \defn{structures} d'espèce $\especeF$ sur $A$, ou encore que ce sont des $\especeF$-structures.
\item[2)]
 La seconde partie de la règle assure que la description de $\especeF[A]$ peut se traduire \defn{naturellement} (au sens précis décrit plus bas) en une description de $\especeF[B]$, chaque fois que $A$ et $B$ ont la même cardinalité. 
Plus précisément, on demande que, pour chaque bijection $\sigma  : A \fleche B$, il y a une bijection
         $$\especeF_{\sigma} :  \especeF[A] \fleche  \especeF[B]$$
décrivant comment transformer les éléments de $\especeF[A]$ en éléments de $\especeF[B]$. On dit que $\especeF_{\sigma}$ est le \defn{transport de structures} le long de $\sigma$.
\end{enumerate} 
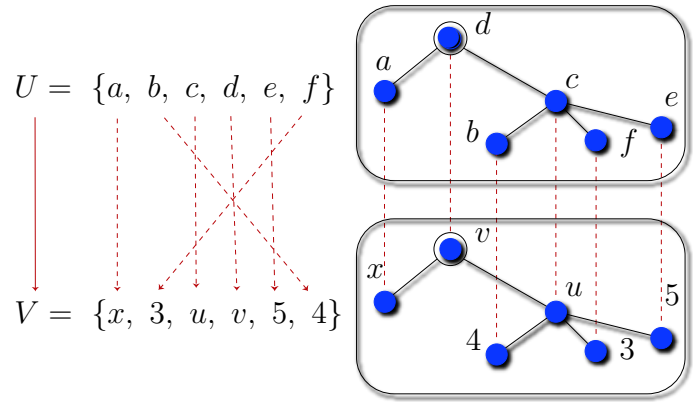
\begin{figure}[ht]
\begin{center}
\scalebox{.6}{\includegraphics{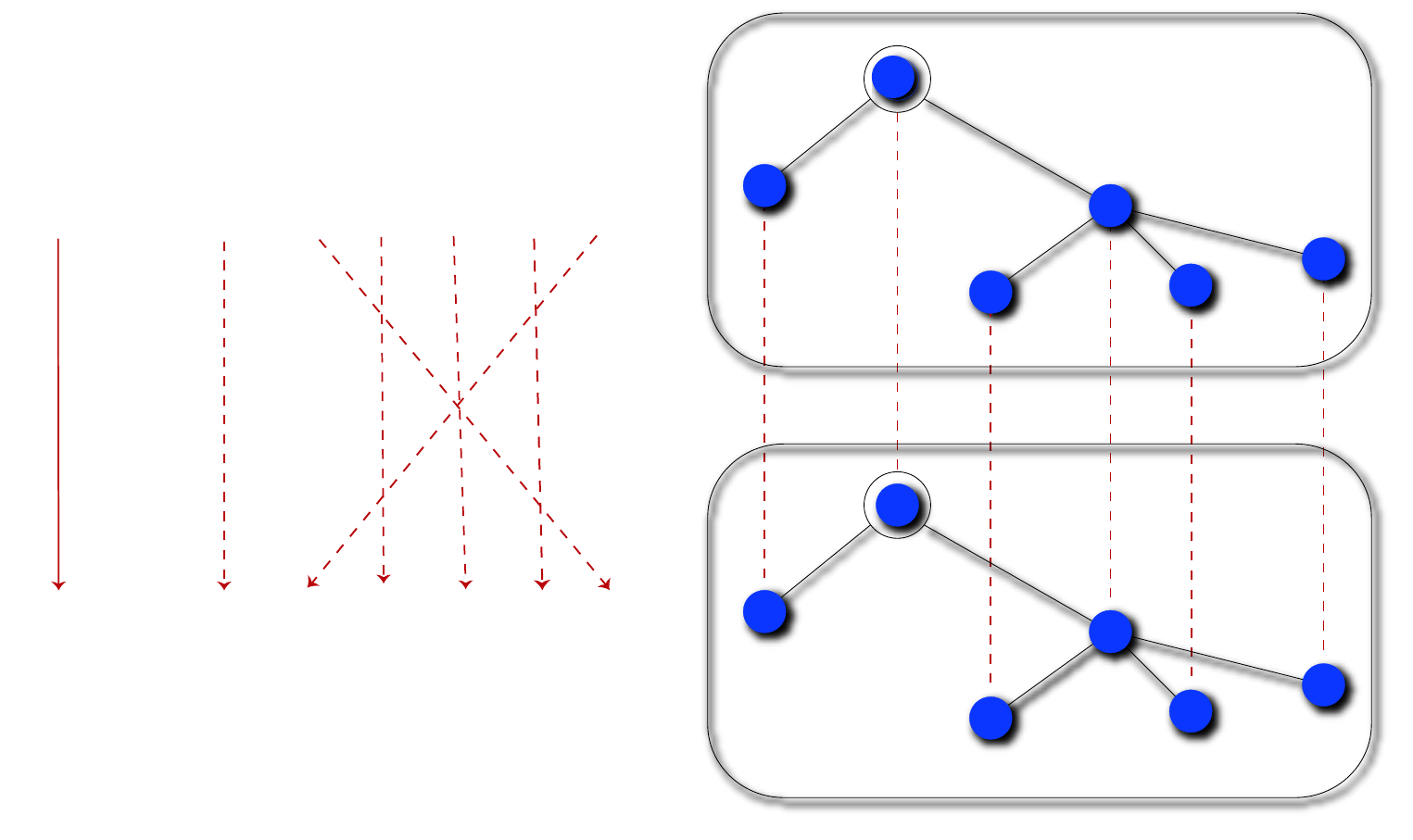}}
 \begin{picture}(0,0)(1.6,.2) 
   \setlength{\unitlength}{1.2mm}
\put(-83,-5){
  \put(52,44){$a$}
  \put(62,36){$b$}
  \put(73,42){$c$}
  \put(63,48){$d$}
  \put(84,40){$e$}
  \put(79,35){$f$}
     \put(51,21){$x$}
  \put(62,13){$4$}
  \put(73,19){$u$}
  \put(63,25){$v$}
  \put(84,18){$5$}
  \put(79,12){$3$}
  \put(12,41){$U=\ \{a,\ b,\ c,\ d,\ e,\ f\}$}
  \put(12,16){$V=\ \{x,\ 3,\ u,\ v,\ 5,\ 4\}$}}
 \end{picture}
 \end{center}\vskip-15pt
\caption{Le ré-étiquetage comme transport.}
\label{fig_transport}
\end{figure}

Très souvent, les éléments de $A$ servent \defn{d'étiquettes} au sein des structures dans $\especeF[A]$, et alors $\especeF_{\sigma}$ correspond à remplacer l'étiquette $x$ par l'étiquette $\sigma(x)$, pour chaque $x$ dans $A$ (voir la Figure~\ref{fig_transport}). Autrement dit, la bijection $\especeF_{\sigma}$ est simplement le \defn{ré-étiquetage} selon $\sigma$. Bien entendu, il y a d'autres possibilités pour $\especeF_{\sigma}$.
Quelque soit le cas, pour exclure les situations potentiellement non désirables (contre intuitives), ainsi que pour donner une définition qui soit techniquement aisément manipulable, on impose que les bijections  $\especeF_{\sigma}$ satisfassent aux conditions de \defn{fonctorialité}:
\begin{enumerate}
\item[i)]  pour chaque ensemble fini $A$, la bijection $\especeF_{\Id_A}$ est l'identité de l'ensemble $\especeF[A]$.
\item[ii)] si $\sigma:A\fleche  B$ et   $\tau:B\fleche  C$, alors  $\especeF_{\tau\circ \sigma}  =   \especeF_{\tau}\circ \especeF_{\sigma}$.
\end{enumerate}

La meilleure façon de s'habituer à cette définition est de donner un certain nombre d'exemples typiques. Nous en verrons toute une liste après la définition suivante.

\subsection*{Série génératrice}
De la définition d'espèce, il suit immédiatement qu'on a forcément 
   $$\card{\especeF[A]}=\card{\especeF[B]},$$
chaque fois qu'il y a une bijection entre les ensembles  $A$ et $B$, puisqu'il y a alors une bijection de $\especeF[A]$ vers $\especeF[B]$.   Le nombre de structures d'espèce $\especeF$  sur un ensemble $A$ est donc simplement une fonction du cardinal de $A$. Pour fin d'énumération, on peut choisir n'importe quel ensemble de cardinal $n$, mais usuellement on prend $A=[n]:=\{1,2,\ \ldots\ ,n\}$. On allège la notation en écrivant simplement $\especeF[n]$ pour l'ensemble des $\especeF$-structures sur l'ensemble $A=[n]$. 

La problématique de la combinatoire énumérative, consiste à calculer la suite de nombres:
$$\card{\especeF[0]},\quad \card{\especeF[1]},\quad \card{\especeF[2]},\quad \card{\especeF[3]},\quad 
 \cdots\quad \card{\especeF[n]},\qquad\cdots$$
Une manière efficace de s'attaquer à cette tâche est de considérer  la \defn{série génératrice} 
 $$\bleu{\especeF(x):=\sum_{n=0}^\infty f_n\,\frac{x^n}{ n!}}\,,$$
où $f_n$ est le nombre d'éléments de $\especeF[n]$.  Un des aspects fascinant de l'approche par séries génératrices, est qu'il y a souvent une fonction classique bien connue dont le développement en série à l'origine correspond précisément à une série génératrice apparaissant ainsi dans un contexte d'énumération. Par calcul direct, on trouve facilement certaines séries.

\section{Exemples d'espèces}
Voici maintenant une liste d'exemples d'espèces et de séries génératrices associées. 
Dans la plupart de ces exemples, le calcul des séries est direct. Cependant certaines séries s'obteniennent par des techniques présentées plus loin.

\subsection*{ Espèces de parties}
On a, par exemple,  l'espèce $\P$ des \defn{parties}, obtenue en posant
\begin{eqnarray*}
   &&\P[A]:=\{U\ |\ U\subseteq A\},\\[4pt]
   &&\P_{\sigma}:\P[A]\fleche \P[B],\qquad \P_{\sigma}(U):=\{ \sigma(x)\ |\ x\in U\}.
\end{eqnarray*}
Le transport d'une partie $U$ de $A$, le long d'une bijection $\sigma:A\fleche B$, est donc l'image de $U$ par $\sigma$. On dénote parfois par $\sigma(U)$ cette image.
On a aussi défini l'espèce des \defn{parties à $k$ éléments}, en posant
  $$\P_k[A]:=\{U\ |\ U\subseteq A,\ \card{U}=k\},$$
avec le transport de structures analogue. Pour la bijection $\sigma(1)=a$, $\sigma(2)=b$, et $\sigma(3)=c$, de $\{1,2,3\}$ vers $\{a,b,c\}$, on a le transport de $\P$-structures:
$$\begin{matrix}
\emptyset & \{1\} & \{2\} & \{3\} & \{1,2\} & \{1,3\} & \{2,3\}& \{1,2,3\} \\
\rouge{\downarrow} & \rouge{\downarrow} &\rouge{\downarrow} &\rouge{\downarrow} &\rouge{\downarrow} &\rouge{\downarrow} &\rouge{\downarrow} &\rouge{\downarrow} \\
\emptyset & \{a\} & \{b\} & \{c\} & \{a,b\} & \{a,c\} & \{b,c\}& \{a,b,c\} 
\end{matrix} $$
Puisqu'il y a $2^n$ structures d'espèce $\P$ sur un ensemble à $n$ éléments, on a la série
   \begin{eqnarray*}
       \P(x)&=&\seriequatre{1}{2}{4}{8}{16}{x}+\ \ldots\ +\bleu{2^n}\,\frac{x^n}{n!}+\ \ldots\\
              &=&\bleu{\exp(2\,x)}.
   \end{eqnarray*}
D'autre part, puisqu'il y a $\binom{n}{k}$ structures d'espèce $\P_k$ sur un ensemble à $n$ éléments, on calcule que
  $$ \P_k(x)=\bleu{\frac{x^k}{k!}\,\exp(x)}.$$

\subsection*{Espèces de graphes}
\begin{figure}[ht]
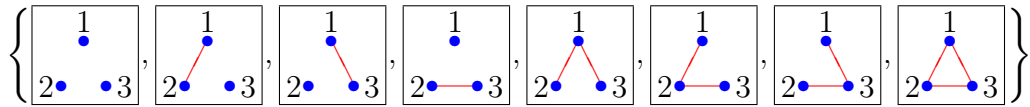

$$
    \left\{\grapheun{1}{2}{3},\graphedeux{1}{2}{3},\graphetrois{1}{2}{3},\graphequatre{1}{2}{3},
   \graphecinq{1}{2}{3},\graphesix{1}{2}{3},\graphesept{1}{2}{3},\graphehuit{1}{2}{3}\right\}
  $$
  \vskip-10pt
\caption{L'ensemble des graphes sur $\{1,2,3\}$.}
\label{les_graphes}
\end{figure}

Considérons l'ensemble $\Gra[A]$ des \defn{graphes simples} sur $A$. C'est là la première partie de la description de l'espèce des graphes simples. On dit des éléments de $\Gra[A]$ que ce sont les structures d'espèce graphe sur $A$, ou simplement les graphes sur $A$. Ainsi, la Figure~\ref{les_graphes} présente les graphes sur $A=\{1,2,3\}$.
\setlength{\graphelength}{.2mm}
Plus techniquement, on pose
 $      \Gra[A]:=\P[\P_2[A]]$.
Reste maintenant à décrire le transport de structures, qui est dans ce cas un ré-étiquetage. 
Pour une bijection $\sigma:A\bijection B$,  le transport $\Gra_{\sigma}$ d'un graphe $G$ dans $\Gra[A]$ s'obtient en posant:
 $$\Gra_{\sigma}(G):=\{\bleu{\{\sigma(a),\sigma(a')\}}\ |\ \bleu{\{a,a'\}}\in G \}$$
Ceci correspond à changer les étiquettes des sommets, en remplaçant $a$ par $\sigma(a)$.
Ainsi, lorsqu'on transporte les graphes sur $\{1,2,3\}$ le long de la  bijection $\sigma:\{1,2,3\}\rightarrow  \{a,b,c\}$, pour laquelle $\sigma(1)=b$, $\sigma(2)=a$ et $\sigma(3)=c$, on obtient le transport:
$$\begin{array}{llll}
    \grapheun{1}{2}{3}\image{\sigma} \grapheun{b}{a}{c},\ & 
    \graphedeux{1}{2}{3}\image{\sigma} \graphedeux{b}{a}{c},&
    \graphetrois{1}{2}{3}\image{\sigma} \graphetrois{b}{a}{c},&
   \graphequatre{1}{2}{3} \image{\sigma} \graphequatre{b}{a}{c},\\[20pt]
    \graphecinq{1}{2}{3}\image{\sigma}  \graphecinq{b}{a}{c},&
    \graphesix{1}{2}{3}\image{\sigma} \graphesix{b}{a}{c},&
    \graphesept{1}{2}{3}\image{\sigma} \graphesept{b}{a}{c},&
    \graphehuit{1}{2}{3} \image{\sigma} \graphehuit{b}{a}{c},\\[20pt]
  \end{array}$$
  
De même, on a l'espèce $\Gro$ des \defn{graphes orientés} définie en posant
  $$\Gro[A]:=\P[A\times A].$$
Le transport de structure correspondant est 
$      \Gro_{\sigma}:=\P_{\sigma\times \sigma}$,
où $\sigma\times \sigma$ est la bijection qui envoie $(a,b)$, dans $A\times A$, sur 
$$\bleu{(\sigma\times \sigma)(a,b):=(\sigma(a),\sigma(b))}.$$
 Autrement dit, l'image par  $\Gro_{\sigma}$ d'un graphe orienté $G$ dans $\Gro[A]$, est le graphe orienté:
   $$\Gro_{\sigma}(G)=\{ (\sigma(a),\sigma(b))\ |\ (a,b)\in G \}.$$
Puisqu'il y a $2^{\binom{n}{ 2}}$ structures d'espèce $\Gra$ sur un ensemble à $n$ éléments, on a la série
   \begin{eqnarray*}
       \Gra(x)
        &=&\seriecinq{1}{}{2}{8}{64}{1024}{x}\\
        &=&\bleu{\sum_{n=0}^\infty 2^{\binom{n}{ 2}}\,\frac{x^n}{ n!} }.
   \end{eqnarray*}
Pour l'espèce des graphes orientés, il y a $2^{(n^2)}$ structures sur un ensemble à $n$ éléments, on a donc la série
   \begin{eqnarray*}
       \Gro(x)&=&\seriecinq{1}{2}{16}{512}{65536}{33554432}{x}\\
        &=&\bleu{\sum_{n=0}^\infty 2^{(n^2)}\,\frac{x^n}{ n!} }.
   \end{eqnarray*}

\subsection*{Espèces d'endofonctions}

 On définit l'espèce $\S$ des \defn{permutations}, en posant
  $$ \S[A]:=\{\tau\ |\ \tau:A\bijection A\}.$$
Le transport de structures, correspondant à une bijection $\sigma:A\fleche B$, se définit en posant
 $\S_{\sigma}(\tau):=\sigma\, \tau\,\sigma^{-1}$,
pour $\tau$ dans $ \S[A]$. Autrement dit, le transport de la permutation $\tau:A\fleche A$ le long de $\sigma:A\fleche B$ donne une permutation $\theta:=\S_{\sigma}(\tau)$ de $B$:
     $$\bleu{\theta(\rouge{\sigma(a)})= \rouge{\sigma(b)},\ssi \tau(\rouge{a})=\rouge{b}}.$$
Plus généralement, pour une endofonction $f:A\fleche A$, on a la conjugaison:
\begin{center}
\begin{picture}(3,4)(0,-.5)
\put(0,3){\bleu{$A$}}  \longuefleche{.8}{3} \put(4,3){\bleu{$A$}}
\flecheverslehaut{0.3}{2.8} \flecheverslebas{4.3}{2.8}
\put(0,0){\bleu{$B$}}   \pointille{0.6}{0} \put(4,0){\bleu{$B$}}
\put(2,3.5){$f$}
\put(-1,1.5){$\sigma^{-1}$}
\put(4.5,1.5){$\sigma$}
\put(1.5,-.7){\bleu{$\sigma f \sigma^{-1}$}}
\end{picture} \end{center}
Cette conjugaison par $\sigma$ donne le transport pour les espèces de fonctions suivantes: 
\begin{enumerate}\itemsep=10pt\itemindent=10pt
  \item[$\bullet$] l'espèce $\Cyc$ des \defn{cycles} (ou des \defn{permutations cycliques}), définie comme
         $$\Cyc[A]:=\{\sigma\ |\ \sigma:A\bijection A,\ \sigma\ {\rm cyclique}\}.$$
  \item[$\bullet$] l'espèce $\Inv$ des \defn{involutions}, définie comme
         $$\Inv[A]:=\{\sigma\ |\ \sigma:A \bijection A,\ \sigma^2=\Id_A\}.$$
  \item[$\bullet$] l'espèce $\Der$ des \defn{dérangements} 
          $$\Der[A]:=\{\sigma\ |\ \sigma:A\bijection A,\ (\forall a)\ \sigma(a)\not=a \},$$
  \item[$\bullet$] l'espèce $\Endo$ des  \defn{endofonctions} 
          $$\Endo[A]:=\{f\ |\ f:A\fleche A\}.$$
\end{enumerate}
Puisqu'il y a $n!$ structures d'espèce $\S$ sur un ensemble à $n$ éléments, on a la série
   \begin{eqnarray*}
       \S(x)&=&\seriequatre{1}{}{2}{6}{24}{x}+\ \ldots\ +\bleu{n!}\,\frac{x^n}{n!}+\ \ldots\\
              &=&\bleu{\frac{1}{1-x}}.
   \end{eqnarray*}
Il y a $(n-1)!$ permutations cycliques sur un ensemble à $n$ éléments. On a donc
   \begin{eqnarray*}
       \Cyc(x)&=&\seriequatre{0}{}{1}{2}{6}{x}+\ \ldots\ +\bleu{(n-1)!}\,\frac{x^n}{n!}+\ \ldots\\
              &=&\bleu{\log\frac{1}{1-x}}.
   \end{eqnarray*}
On montrera plus loin que le série génératrice de l'espèce des involutions est 
   \begin{eqnarray*}
       \Inv(x)&=&\seriecinq{1}{}{2}{4}{10}{26}{x}\\
              &=&\bleu{\exp(x+x^2/2)},
   \end{eqnarray*}
et que celle des dérangements est
   \begin{eqnarray*}
       \Der(x)&=&\seriecinq{1}{0}{}{2}{9}{44}{x}\\
              &=&\bleu{\frac{\exp(-x)}{1-x}}.
   \end{eqnarray*}
Comme il y a $n^n$ endofonctions sur un ensemble à $n$ éléments, on a
   \begin{eqnarray*}
       \Endo(x)&=&\seriecinq{1}{}{4}{27}{256}{3125}{x}\\
        &=&\bleu{\sum_{n=0}^\infty n^n\,\frac{x^n}{ n!} },
   \end{eqnarray*}
 
\subsection*{L'espèce des listes}

 On définit encore l'espèce $\Liste$ des \defn{listes}  (ou \defn{ordres linéaires}):
  $$ \Liste[A]:=\{ \ell\ |\ \ell:[n]\ \bijection\ A\}.$$
Autrement dit, une liste des éléments de $A$ s'identifie à une fonction 
     $$\bleu{\ell=\{ (1,x_1), (2,x_2), (3,x_3),\ \ldots\, (n,x_n) \}}.$$
On peut alors définir le transport de structures le long de $\sigma$ tout simplement comme
  $$\Liste_{\sigma}(\ell):=\sigma\circ \ell.$$
Puisqu'il y a $n!$ structures d'espèce $\Liste$ sur un ensemble à $n$ éléments, on a la série
   \begin{eqnarray*}
       \Liste(x)&=&\seriequatre{1}{}{2}{6}{24}{x}+\ \ldots\ +\bleu{n!}\,\frac{x^n}{n!}+\ \ldots\\
              &=&\bleu{\frac{1}{1-x}}.
   \end{eqnarray*}

\subsection*{L'espèce des partitions}
Rappelons qu'une partition d'un ensemble $A$ est un ensemble $P$ de parties non vides de $A$, deux-à-deux disjointes, dont la réunion est $A$. On a  l'espèce $\Part$ des \defn{partitions} en prenant
   $$\Part[A]:=\{ P\ |\ P\ \hbox{est une partition de}\ A\},$$
   et en considérant de pair le transport de structures défini par
  $$ \Part_{\sigma}(P):=\{ \sigma(B)\ |\ B\in P\}$$
pour $P$ dans $\Part[A]$. Comme on va le voir plus loin, l'espèce des partitions joue un rôle crucial dans la définition de l'opération de substitution d'espèces.
On a la série (voir plus loin pour la justification)
   \begin{eqnarray*}
       \Part(x)&=&\seriecinq{1}{}{2}{5}{15}{52}{x}\\
              &=&\bleu{\exp(\exp(x)-1)},
   \end{eqnarray*}
dont les coefficients sont les nombres de Bell.


\subsection*{Espèces auxiliaires simples.}\label{def_autres_especes}
Voici quelques autres exemples, parfois presque trop simples, qui seront utiles pour la suite. Dans chaque cas on décrit l'effet de l'espèce sur un ensemble $A$, et le transport pour $\sigma:A\bijection B$. Cependant, lorsque le transport de structures est assez évident, il n'est pas nécessairement décrit explicitement.

\begin{enumerate}\itemsep=10pt\itemindent=10pt
  \item[(1)] $\Zero[A]=\emptyset$, \quad  c'est l'espèce dite \defn{vide}. Le transport de structure est l'unique bijection de l'ensemble vide vers l'ensemble vide. Comme il n'y a aucune structures d'espèce $\Zero$ sur un ensemble de cardinal $n$, la série génératrice associée est
    \begin{eqnarray*}
        \bleu{\Zero(x)}&=& \seriecinq{0}{0}{0}{0}{0}{0}{x}\\
                               &=&\bleu{0}.
     \end{eqnarray*}
  \item[(2)] $\Ens[A]=\{A\}$,\quad c'est  l'espèce des \defn{ensembles} (avec une seule structure\footnote{Attention ici à la différence entre l'ensemble $A$, qui peut contenir plusieurs éléments, et l'ensemble $\{A\}$ qui ne contient qu'un seul élément, à savoir l'ensemble $A$.} sur tout ensemble fini). On pose $\Ens_{\sigma}(A):=B$.  Il y a exactement une structure d'espèce $\Ens$ sur chaque ensemble de cardinal $n$, on a donc
   \begin{eqnarray*}
       \Ens(x)&=& 1+ x+ \frac{x^2}{2}+ \frac{x^3}{6}+ \frac{x^4}{24}+\ \ldots\ + \frac{x^n}{n!}+\ \ldots\\
              &=&\bleu{\exp(x)}.
   \end{eqnarray*}

   \item[(3)]  $  \Un[A] =\begin{cases}
      \{A\} & \text{si } \ \card{A}=0,\\[4pt]
     \emptyset & \text{sinon}.
\end{cases} $ \\
C'est l'espèce \defn{caractéristique de l'ensemble vide}. Le seul ensemble qui admet une structure d'espèce $\Un$ est l'ensemble vide, on a donc
    \begin{eqnarray*}
        \bleu{\Un(x)}&=& \seriecinq{1}{0}{0}{0}{0}{0}{x}\\
                               &=&\bleu{1}.
     \end{eqnarray*} 
  \item[(4)]  $  \SingletonX[A] = \begin{cases}
   \{A\}    & \text{si}\  \card{A}=1, \\[4pt]
   \emptyset   & \text{sinon}.
\end{cases} $ \\
C'est  l'espèce \defn{caractéristique des singletons}. Les seuls ensembles qui admettent des structures d'espèce $\SingletonX$ sont les singletons, et il y a alors une seule telle structure. On a donc
     \begin{eqnarray*}
        \bleu{\SingletonX(x)}&=& \seriecinq{0}{}{0}{0}{0}{0}{x}\\
                               &=&\bleu{x}.
     \end{eqnarray*}
  \item[(5)] Pour toute espèce $\especeF$, on peut considérer l'espèce $\especeF\nonvide$ qui est sa restriction aux ensembles non-vides:
$$\especeF\nonvide[A] =\begin{cases}
      \especeF[A] & \text{si } \ \card{A}\not= 0,\\[4pt]
     \emptyset & \text{sinon}.
\end{cases} $$
\end{enumerate}

\section{Isomorphisme/égalité d'espèces.}
L'isomorphisme  d'espèces est la formulation précise de la notion de \og bijection naturelle\fg\ entre structures combinatoires. Informellement, une bijection est dite naturelle si sa description ne dépend pas de propriétés particulières des éléments (grandeur, valeur, ...) qui interviennent dans les structures mises en bijection. Le critère technique plus explicite consiste à demander que la bijection soit \og compatible\fg\ avec les transports de structures. Le lecteur intéressé pourra consulter la monographie~\cite{BLL} pour une présentation plus détaillée. 

Cependant, la notion d'isomorphisme est essentielle pour les développements les plus significatifs de la théorie des espèces, par exemple pour l'énumération de structures non étiquettées, ou encore pour établir un lien clair avec la théorie de la représentation de groupes des permutations. En fait, c'est cet aspect de la théorie qui la rend vraiment originale par rapport au autres théories qui ont été proposées pour rendre compte de l'approche par séries génératrices à la combinatoire énumérative.

Il est souvent agréable d'écrire $\especeF= \especeG$, pour signifier qu'il existe un isomorphisme entre les espèces $\especeF$ et $\especeG$ sans avoir nécessairement à spécifier lequel. Cela est certainement un abus de langage, mais le contexte rend toujours justice aux exigences de rigueur. Du point de vue tout-à-fait pratique, on établit qu'il y a égalité entre deux espèces en décrivant comment transformer les structures de la première espèce en structures de la seconde, d'une façon clairement inversible. Très souvent cette transformation est facile à comprendre en l'illustrant judicieusement par un exemple assez générique.

\section{Somme et produit.}\label{sec_som_prod}
Nous allons maintenant développer l'un des aspects les plus intéressant de la théorie des espèces: \og l'algèbre des espèces\fg. Plus précisément, nous allons introduire des opérations entre espèces qui permettent de manipuler, de construire et de décomposer les espèces. Nous allons constater que plusieurs arguments combinatoires se formulent agréablement en terme de ces opérations. 

On a d'abord la \defn{somme} d'espèces:
 $$(\especeF+\especeG)[A]:=\especeF[A]+\especeG[A],\qquad \hbox{(somme disjointe)}$$
Autrement dit, une structure d'espèce $\especeF+\especeG$ sur $A$ est soit une structure d'espèce $\especeF$, soit une structure d'espèce $\especeG$. Pour $\sigma:A\bijection  B$, on pose
$$   (\especeF+\especeG)_{\sigma}(t):=\begin{cases}
    \especeF_{\sigma}(t)  & \text{si}\quad t\in \especeF[A], \\
    \especeG_{\sigma}(t)  & \text{si}\quad  t\in \especeG[A].\end{cases}$$
On vérifie facilement que la somme d'espèce est associative et commutative (à isomorphisme d'espèce près), et que l'espèce $\Zero$ agit comme élément neutre. 

Plus intéressant encore est le \defn{produit}, $\especeF\cdot \especeG$, d'espèces qui est défini sur un ensemble $A$ comme
 $$(\especeF\cdot \especeG)[A]:=\sum_{A=B+C} \especeF[B]\times \especeG[C],$$
où la sommation correspond à l'union disjointe.
  Le transport de structures d'espèce $\especeF\cdot \especeG$ se définit en posant
 $$(\especeF\cdot \especeG)_{\sigma}(f,g)=(\especeF_{\sigma}(f),\especeG_{\sigma}(g)).$$
Nous allons voir des exemples un peu plus loin.
\section{Dessins de structures génériques.}\label{sec_dessinsesp}
Pour ne pas tomber dans un excès de formalisme risquant d'obscurcir inutilement les raisonnements, il est agréable et efficace de décrire les constructions en terme de dessins qui présentent un exemple typique de structure de l'espèce considérée. C'est d'ailleurs ce que nous avons fait pour l'espèce des graphes orientés ($\Gro$). 

Il est intéressant d'étendre cette habitude aux nouvelles espèces introduites par des opérations combinatoires.
Pour systématiser un peu cette façon de faire, et 
afin de décrire diverses opérations et manipulations générales sur les espèces de façon plus conviviale, 
on adopte souvent une présentation \og par le dessin\fg\ des
structures \og génériques\fg\ d'espèces abstraites $\especeF$,   $\especeG$, etc. 

\begin{figure}[h!]
  \begin{center}
  \scalebox{.8}{\includegraphics{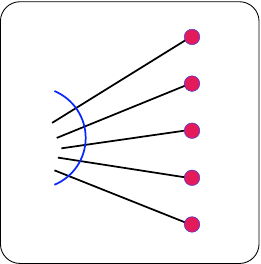}}
  \begin{picture}(0,0)(0,0)
    \put(-4,.5){\large$\especeF$}
  \end{picture}
\qquad 
 \scalebox{.4}{\includegraphics{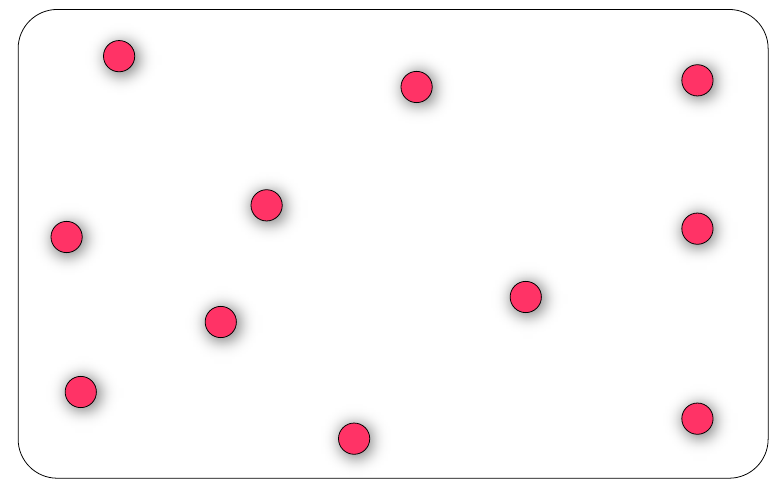}}
  \begin{picture}(0,0)(0,0)
    \put(-4,1.5){\huge$\especeF$}
  \end{picture}
   \end{center}
   \vskip-10pt
   \caption{Représentations possibles de structure générique d'espèce $\especeF$.}\label{figgeneric}
\end{figure}
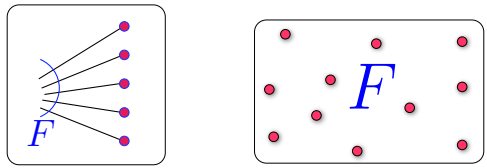

Le premier dessin de la Figure~\ref{figgeneric} représente un élément typique de l'ensemble $\especeF[A]$. Ici les divers éléments de $A$  apparaissent comme des tiges auxquelles sont rattachées des points rouges. L'arc de cercle (avec son indice $\especeF$) symbolise abstraitement le fait qu'on a une certaine structure d'espèce $\especeF$ sur ces éléments.
Une autre possibilité de représentation consiste  simplement à superposer le symbole $\especeF$ sur l'ensemble sous-jacent, tel qu'illustré aussi à la figure~\ref{figgeneric}.

\begin{figure}[h!]
  \begin{center}
\scalebox{.5}{\includegraphics{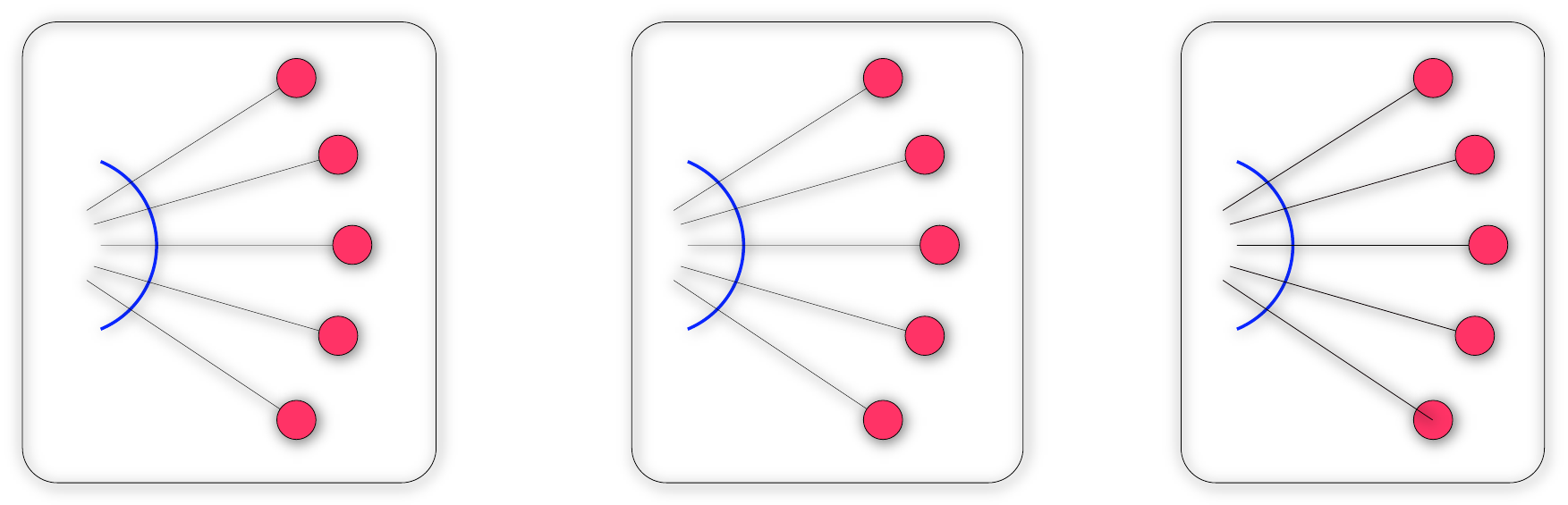}}
  \begin{picture}(0,0)(0,0)
    \put(-17.5,1.3){\small$\especeF+\especeG$}
    \put(-12.5,2.5){\Large$=$}
    \put(-10,1){\small $\especeF$}
    \put(-6,2.5){\large \bleu{ou}}
        \put(-4,1){\small $\especeG$}
      \end{picture}
   \end{center}
   \vskip-10pt
   \caption{Une structure d'espèce $\especeF + \especeG$.}\label{figaddition}
\end{figure}

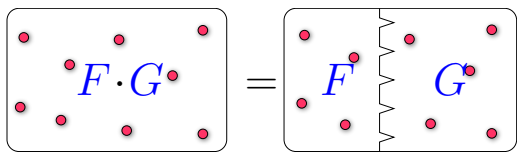
\begin{figure}[h!]
  \begin{center}
\scalebox{.4}{\includegraphics{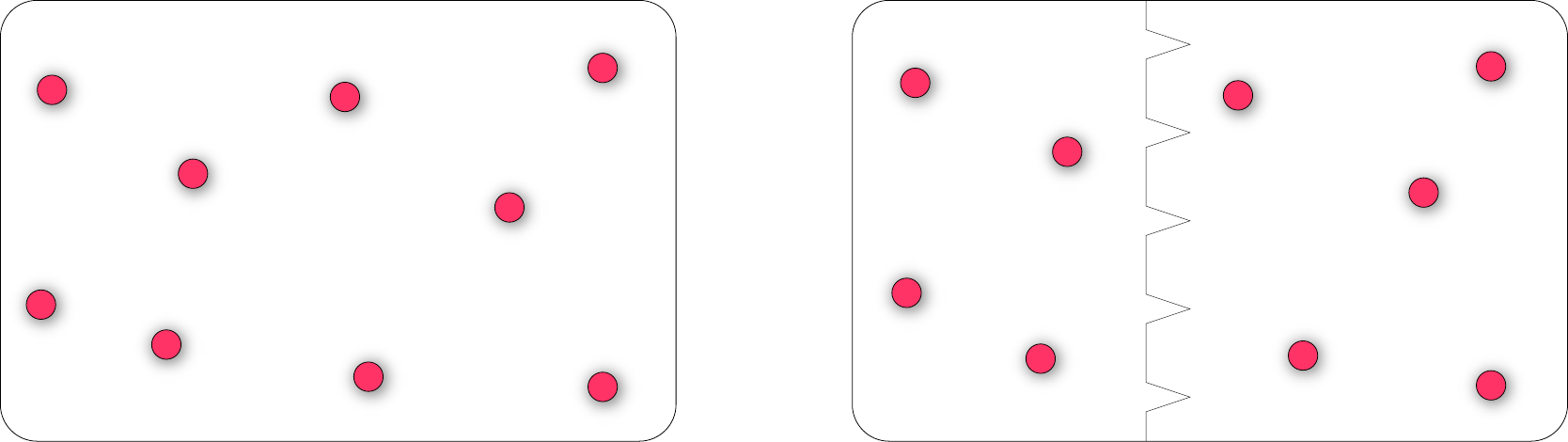}}
  \begin{picture}(0,0)(0,0)
    \put(-12,1.5){\Large$\especeF\!\cdot\! \especeG$}
    \put(-7.5,1.5){\Large$=$}
    \put(-5.5,1.5){\Large$\especeF$}
    \put(-2.5,1.5){\Large$\especeG$}
      \end{picture}
   \end{center}
   \vskip-10pt
   \caption{Une première façon de présenter une structure d'espèce $\especeF\cdot \especeG$.}\label{figprodun}
\end{figure}

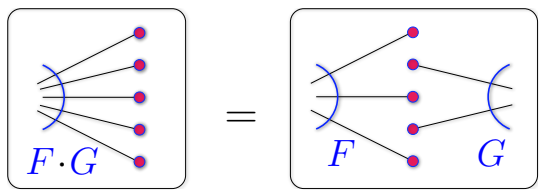
\begin{figure}[h!]
  \begin{center}
\scalebox{.3}{\includegraphics{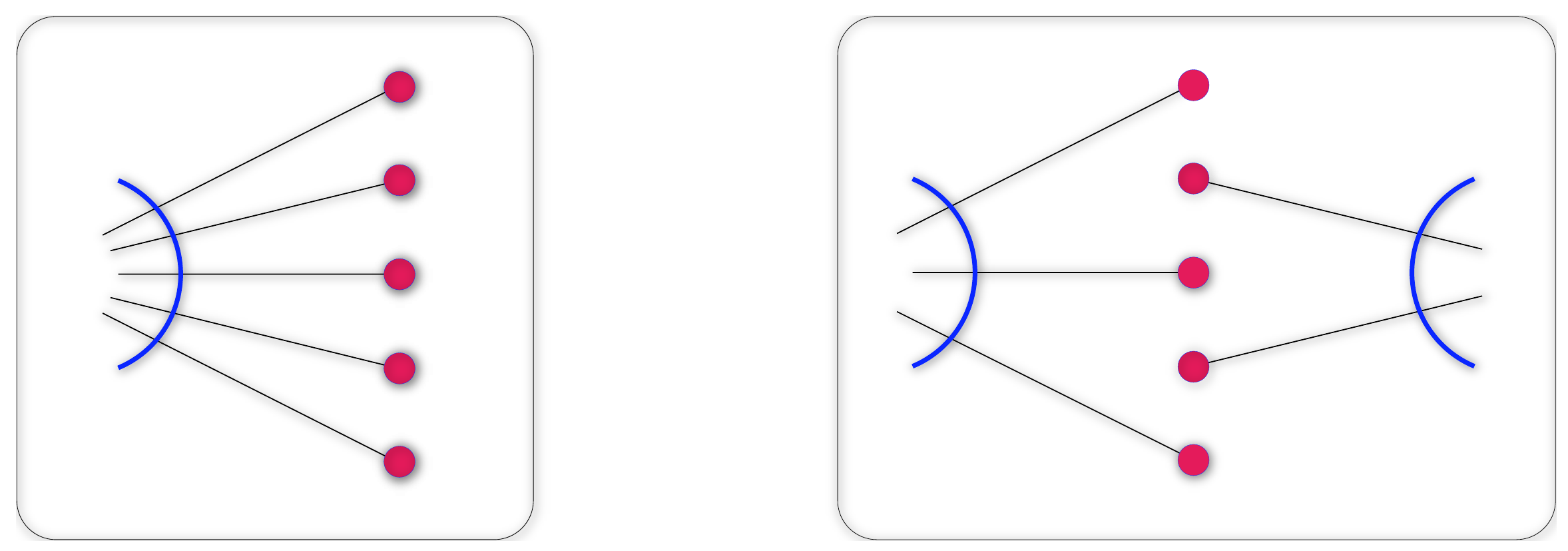}}
  \begin{picture}(0,0)(0,0)
 \put(0,-.3){   \put(-14,.8){\large$\especeF\!\cdot\! \especeG$}
    \put(-8.7,2){\Large$=$}
    \put(-6,1){\large$\especeF$}
    \put(-2,1){\large$\especeG$}}
  \end{picture}
   \end{center}
   \vskip-10pt
   \caption{Une seconde façon de présenter une structure d'espèce $\especeF\cdot \especeG$.}\label{figproddeux}
\end{figure}

On peut alors reformuler, en ces termes plus imagés, la définition d'addition d'espèces, comme on l'a fait à la figure~\ref{figaddition}. 
De même, on peut reformuler de façon imagée la définition de produit d'espèces  $\especeF\cdot \especeG$, soit par la figure~\ref{figprodun}, soit par la figure~\ref{figproddeux}. 
La seconde présentation a l'avantage de mieux mettre en évidence qu'on doit considérer tous les découpages ($A=B+C$) de l'ensemble $A$, sans avoir à \og déplacer\fg\ les 
éléments de $A$.

Il y a de nombreux avantages à présenter les choses de ce point de vue plus imagé. En un certain sens, cela est tout aussi
rigoureux si on envisage cette utilisation de figures comme un nouveau genre de formalisme bi-dimensionnel. En cas de doute, il est toujours possible de revenir à une description rigoureusement ensembliste par une traduction presque directe.

\section{Passage aux séries.}\label{sec_sergenesp}
 Observons qu'on a immédiatement $\especeF(x)=\especeG(x)$, si $\especeF$ et $\especeG$ sont des espèces isomorphes. Nous allons maintenant expliquer \og pourquoi\fg\  la méthode des séries génératrices est si puissante. En fait, en plus d'être compatible avec la notion d'égalité (isomorphisme), une des caractéristiques du passage à la série génératrice est d'être compatible avec les opérations entre espèces. Plus précisément,
 \begin{thm}\label{thm_fondamental}
 \bleu{Pour toutes espèces $\especeF$ et $\especeG$, on a}
\begin{eqnarray*}
    (\especeF+\especeG)(x)&=&\especeF(x)+\especeG(x),\\
   (\especeF\cdot \especeG)(x)&=&\especeF(x)\,\especeG(x),
    \end{eqnarray*}
 \end{thm}
La démonstration de ces égalités est directe (modulo les définitions).

\section{Exemples de calcul de séries.}\label{sec_exemples}
Un exemple  intéressant concerne l'espèce $\Der$ des \defn{dérangements}. Toute permutation $\sigma$ de $A$ se décompose (voir la Figure~\ref{figderang}) en: 
\begin{enumerate}
 \item[i)] l'ensemble, disons $B$, de ses points fixes, et 
 \item[ii)] un dérangement de $A\setminus B$.
\end{enumerate} 
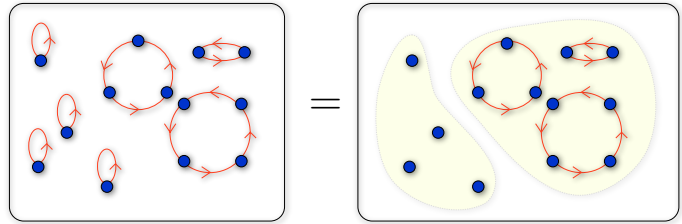
\begin{figure}[ht]
  \begin{center}
  \scalebox{.5}{\includegraphics{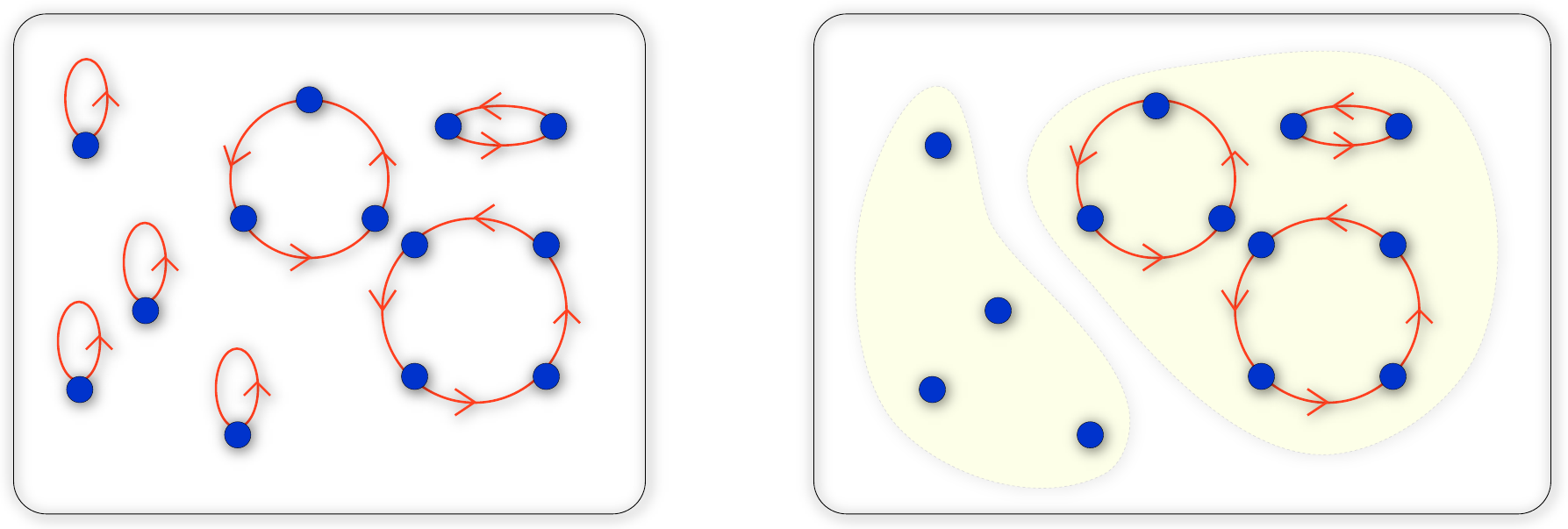}}
    \begin{picture}(0,0)(0,0)
             \put(-10.3,3){\Large$=$}
   \end{picture}
   \end{center}
   \vskip-10pt
   \caption{Une permutation $=$ points fixes et dérangement.}\label{figderang}
\end{figure}
Ceci se traduit en un isomorphisme d'espèces: 
$  \S =  \Ens\cdot \Der$, qu'on peut lire:

\begin{center}
\og \bleu{Une permutation est constituée\\ d'un ensemble de points fixes accompagné d'un dérangement.}\fg
\end{center}

 On en déduit l'identité de séries $\S(x)=\Ens(x)\,\Der(x)$. Comme on connaît les séries génératrices pour les espèces $\S$ et $\Ens$, on trouve
 $$ \bleu{\Der(x)=\frac{\exp(-x)}{ 1-x}}.$$
Il en découle une démonstration simple  de la formule souvent citée 
  $$\bleu{d_n=n!\,\sum_{k=0}^n \frac{(-1)^k}{ k!}},$$
pour le nombre de dérangements de $n$ éléments.
On peut même définir implicitement une espèce comme solution d'une équation. Par exemple, l'espèce $\Bin$ des \defn{arbres binaires} (étiquetés) est l'{unique}\footnote{Il y a un analogue (à la mode des espèces) du théorème des fonctions implicites qui assure l'existence et l'unicité de cette solution.} solution de l'équation 
\begin{equation}\label{abc}
  \Bin=  1+\SingletonX\cdot  \Bin^2.
\end{equation}
Cette équation  admet la présentation imagée:
\begin{center}
   \scalebox{.8}{\includegraphics{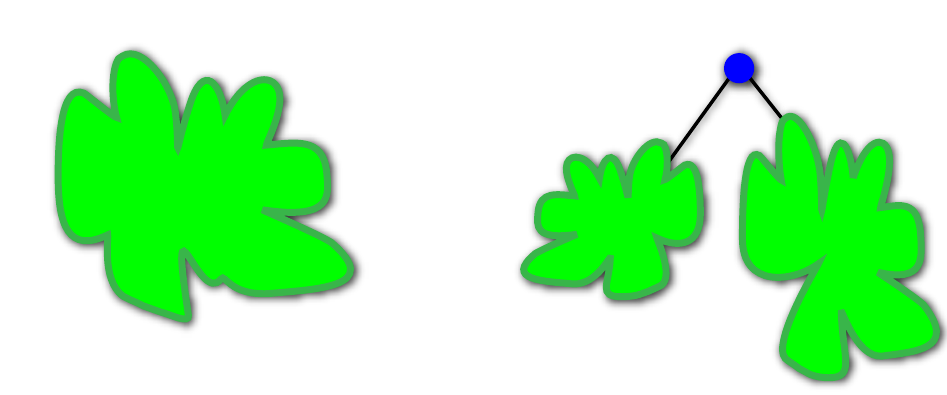}}
   \begin{picture}(0,0)(0,0)
      \put(-13,3){\Large$\Bin$}
       \put(-6,2.3){\Large$\Bin$}
       \put(-2,2.3){\Large$\Bin$}
             \put(-9,3.5){\Large$=$}
      \put(-4,6){\Large$x$}
   \end{picture}
 \end{center}
 Autrement dit, 
   
   \begin{center}
     \og \bleu{Un arbre binaire (non vide) est constitué d'une racine à laquelle sont attachés deux arbres binaires, un à gauche et un à droite.}\fg 
   \end{center}
   
Il découle de l'équation (\ref{abc}) que
   $$\Bin(x)=1+x\,\Bin(x)^2$$
dont la solution  (la seule admettant un développement en séries à l'origine) est
\begin{equation}
     \bleu{\Bin(x)=\frac{1-\sqrt{1-4\,x}}{ 2\,x}=\sum_{n\geq 0} n!\,C_n\,\frac{x^n}{n!}}.
 \end{equation}
Il s'ensuit que le nombre d'arbres binaires (étiquetés) sur $[n]$ est $n!\,C_n$, où $C_n$ est le $n^{\rm e}$ nombre de Catalan $\frac{1}{ n+1}\binom{2\,n}{ n}$.  

\section{Substitution.}\label{sec_substesp}
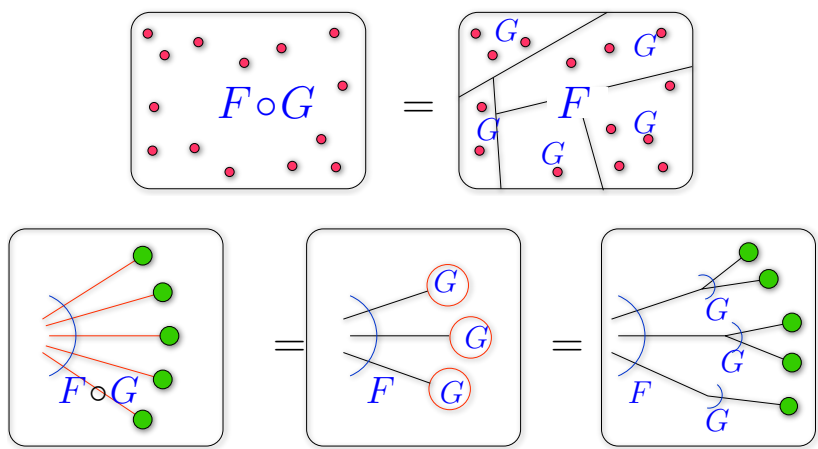
\begin{figure}[ht]\begin{center}
\scalebox{.5}{\includegraphics{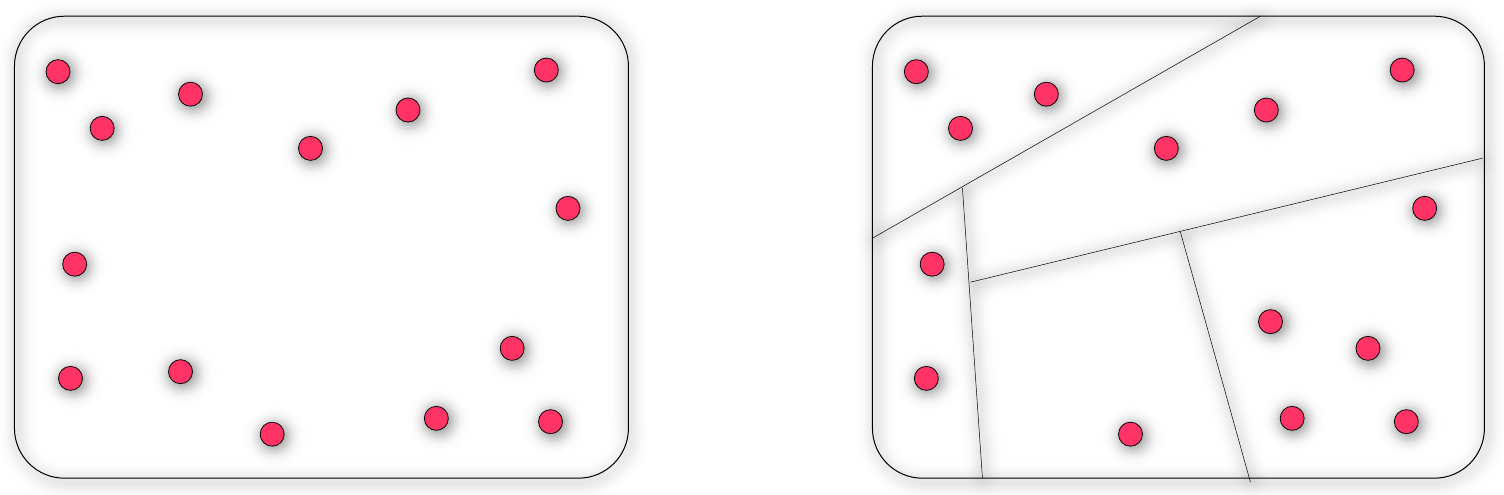}}
 \begin{picture}(0,0)(0,0)\setlength{\unitlength}{4.1mm}
    \put(-16,2.5){\Large\bleu{$\especeF\!\circ\! \especeG$}}
    \put(-10,2.5){\Large$=$}
    \put(-4.3,2.5){\linethickness{8mm}\blanc{\line(0,1){1}}}
    \put(-5,2.5){\Large $\especeF$}
    \put(-2.5,2){$\especeG$}
    \put(-2.5,4.5){$\especeG$}
    \put(-7.6,1.8){$\especeG$}
    \put(-5.5,1){$\especeG$}
         \put(-7,5){$\especeG$}
          \end{picture}
\vskip10pt
\scalebox{.5}{\includegraphics{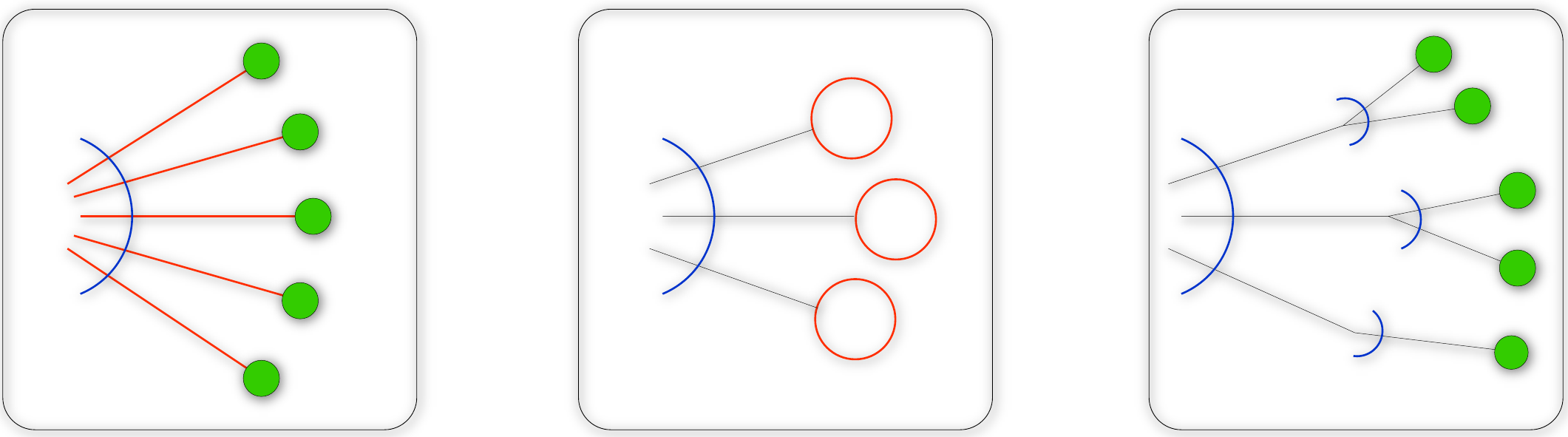}}
 \begin{picture}(0,0)(0,0)\setlength{\unitlength}{4.1mm}
    \put(-25,1.5){\large$\especeF\!\circ\! \especeG$}
    \put(-18,3){\Large$=$}
    \put(-9,3){\Large$=$}
    \put(-15,1.5){\large$\especeF$}
    \put(-12.6,1.5){$\especeG$}
    \put(-11.8,3.3){$\especeG$}
    \put(-12.8,5.1){$\especeG$}
   \put(-6.5,1.5){$\especeF$}
    \put(-4,.6){$\especeG$}
    \put(-3.5,2.5){$\especeG$}
    \put(-4,4.2){$\especeG$}
      \end{picture}
\end{center}
\vskip-10pt
\caption{Diverses représentations de la substitution.}\label{fig_subs1}
\end{figure}

L'opération qui est certainement la plus intéressante est celle de \defn{substitution}.  Elle a lieu entre un espèce  $\especeG$ telle que\footnote{Nous allons voir plus loin pourquoi cette condition est nécessaire lorsqu'on veut définir la substitution en toute généralité.} $\especeG[\emptyset]=\emptyset$, et une espèce $\especeF$ quelconque. On pose
\begin{equation}
   (\especeF\circ \especeG)[A]:=\sum_{\pi\in\Part[A]} \especeF[\pi]\times\prod_{B\in\pi} \especeG[B],
\end{equation}
où $\Part[A]$ désigne l'ensemble des partitions de $A$. On écrit aussi $\especeF(\especeG)$ pour $\especeF\circ \especeG$. On peut représenter graphiquement une structure typique d'espèce $\especeF\circ \especeG$ par l'un des dessins de la Figure~\ref{fig_subs1}.

Il est amusant de constater que nos notations ont été judicieusement choisies pour faire en sorte que 
$ \especeF =  \especeF(\SingletonX)$.
Bien entendu, comme pour les autres opérations, on vérifie par calcul direct qu'il y a compatibilité de la substitution avec le passage aux séries génératrices:
\begin{thm}
   \bleu{Pour tout espèce $\especeF$, et tout espèce $\especeG$ telle que $\especeG[\emptyset]=\emptyset$, on a}
   \begin{equation}
          (\especeF\circ \especeG)(x)=\especeF(\especeG(x)).
   \end{equation}
\end{thm}
 Ainsi si
    $$\especeF(x)=\sum_{n=0}^\infty f_n\,\frac{x^n}{ n!},\qquad {\rm et}\qquad 
         \especeG(x)=\sum_{k=1}^\infty g_k\,\frac{x^k}{ k!}  $$
alors
\begin{eqnarray*}
  \especeF(\especeG(x))  & = & \sum_{n=0}^\infty f_n\,\frac{\especeG(x)^n}{ n!}.\\
               & = & \sum_{n=0}^\infty \frac{f_n}{ n!}\,\left(\sum_{k=1}^\infty g_k\,\frac{x^k}{ k!}\right)^n\\
               & = & f_0+f_1g_1\,x+(f_2g_1^{2}+f_1g_2)\frac{{x}^{2}}{ 2!}
                  +(f_3g_1^{3}+3\,f_2g_1g_2+f_1g_3)\frac{{x}^{3}}{ 3!}\\
              & & \qquad    +  (f_4g_1^{4}+6\,f_3g_1^{2}g_2+4\,f_2g_1g_3+3\,f_2{g_2}^{2}+f_1g_4 )\frac{{x}^{4}}{ 4!}
          + \ldots
\end{eqnarray*}
On interprète combinatoirement le terme $6\,f_3\,g_1^2\,g_2$, apparaissant dans le coefficient de $x^4/4!$ du développement de $\especeF(\especeG(x))$, de la façon suivante. Parmi les $15$ partitions d'un ensemble à $4$ éléments (ici $\{a,b,c,d\}$), il y en a $6$ avec deux parties à $1$ élément, et une à $2$ éléments:
\begin{eqnarray*}
    \{\{a\},\{b\},\{c,d\}\},\ &\{\{a\},\{c\},\{b,d\}\},\ &\{\{a\},\{d\},\{b,c\}\},\\
    \{\{b\},\{c\},\{a,d\}\},\ &\{\{b\},\{d\},\{a,c\}\},\ &\{\{c\},\{d\},\{a,b\}\}
\end{eqnarray*}
Le nombre de structures d'espèce $(\especeF\circ \especeG)$ pour chacune de ces partitions est $f_3\,g_1^2\,g_2$. On choisit, en effet, une structure d'espèce $\especeG$ sur chacune des parties ce qui peut se faire de $g_1^2\,g_2$ façons. Puis on choisit une structure d'espèce $\especeF$ sur l'ensemble à $3$ éléments dont les éléments sont ces parties, ce qui peut se faire de $f_3$ façons.

\subsection*{Exemples de substitutions.}
\begin{figure}[ht]
\begin{center}
\scalebox{.5}{\includegraphics{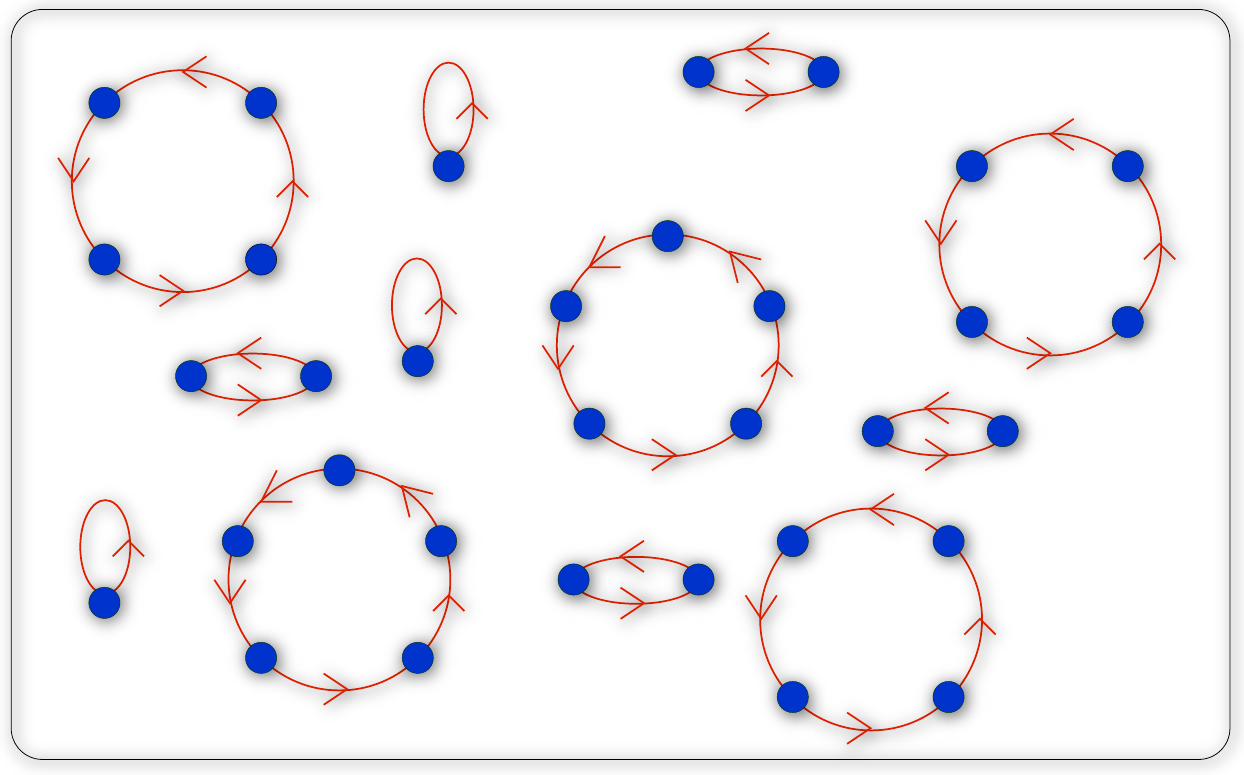}}
\end{center}
\vskip-10pt
\caption{Une permutation est un ensemble de cycles.}\label{ens_de_cycles}
\end{figure}

Puisqu'une partition est un \og ensemble d'ensembles non-vides\fg, on peut décrire 
l'espèce des \defn{partitions} par l'égalité
  $$\Part=
   (\Ens\circ \Ens^+)$$
où $\Ens^+$ désigne l'espèce des ensembles non vides.
On obtient alors
\begin{eqnarray*}
   \Part(x)& = & \exp(\exp(x)-1)\\
                       & = & 1+ x+2\,{\frac {{x}^{2}}{2!}}+5\,{\frac {{x}^{3}}{3!}}+15\,{\frac {{x}^{4}}{4!}}+52\,{
\frac {{x}^{5}}{5!}}+203\,{\frac {{x}^{6}}{6!}}+877\,{\frac {{x}^{7}}{7!}}+4140\,{\frac {
{x}^{8}}{8!}}+\ldots 
\end{eqnarray*}
Les $5$ structures d'espèce $\Part$ sur l'ensemble  $\{a,b,c\}$ sont:
   $$\{\{a,b,c\}\},\ \{\{a,b\},\{c\}\},\ \{\{a,c\},\{b\}\},\ \{\{a\},\{b,c\}\},\ \{\{a\},\{b\},\{c\}\}$$
La Figure~\ref{ens_de_cycles} montre \og comment\fg\  l'espèce des permutations est égale à l'espèce des \defn{ensembles de cycles} $(\Ens\circ \Cyc )$
où $\Cyc $ est l'espèce des \defn{cycles} (\defn{permutations circulaires}). On a donc
$$ e^{\Cyc(x)} = \S(x),$$
d'où
    $$\Cyc(x)=\log\frac{1}{ 1-x},$$
    comme on l'a déjà annoncé.
Ce résultat est consistant avec le fait que $\card{\Cyc[n]}=(n-1)!$, $n\geq 1$, puisqu'on a l'égalité
   $$-\log(1-x)=\sum_{n\geq 1}^\infty \frac{x^n}{ n}.$$

\section{L'espèce des arborescences.}\label{sec_especearbo}
\begin{figure}[ht]
\begin{center}
\scalebox{.6}{\includegraphics{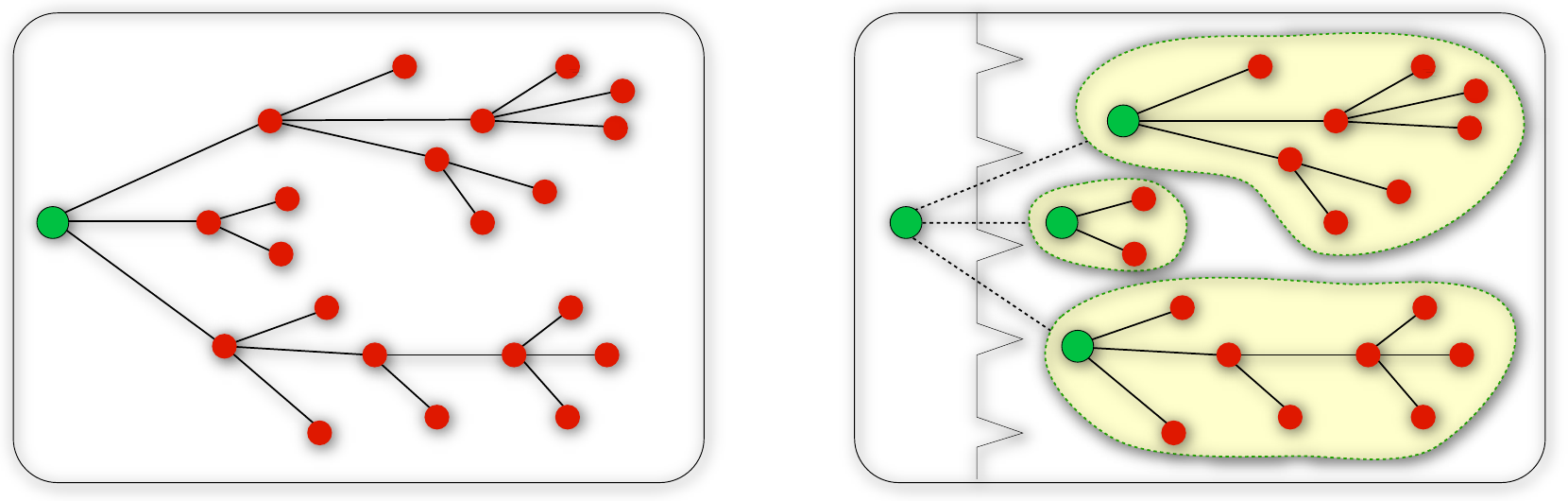}}
\begin{picture}(0,0)(0,0)
  \put(-11,3){\Large{$=$}}
\end{picture}
\end{center}
\caption{Une arborescence \og est\fg\ un ensemble d'arborescences attachées à une racine.}\label{arbor}
\end{figure}
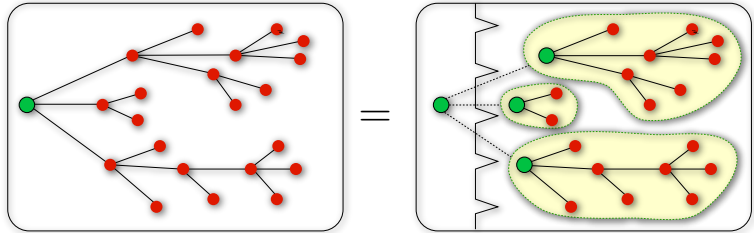

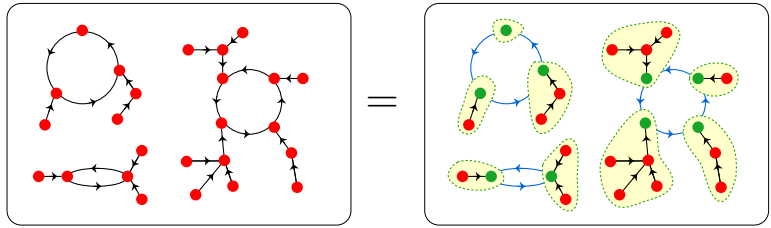
\begin{figure}[ht]
\begin{center}
\scalebox{.8}{\includegraphics{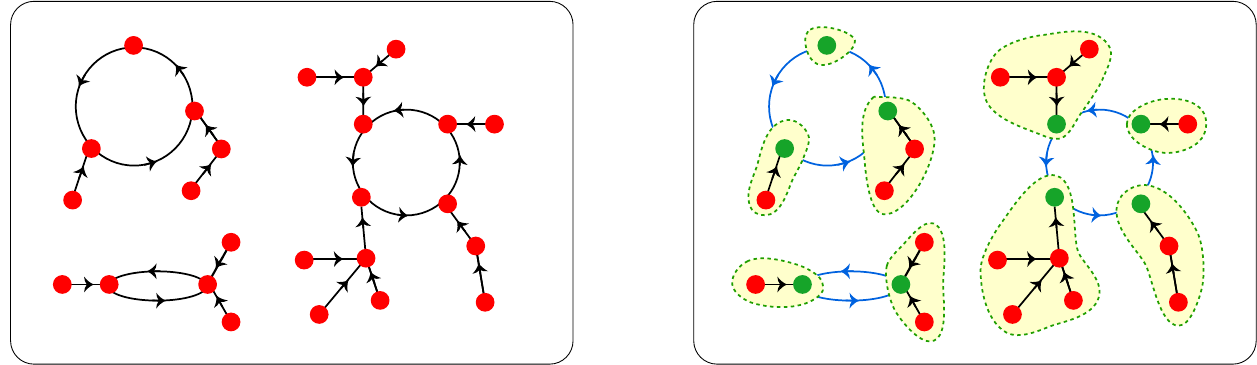}}
\begin{picture}(0,0)(0,0)
  \put(-11,3){\Large{$=$}}
\end{picture}
\end{center}
\vskip-10pt
\caption{La décomposition d'une endofonction}
\label{endo_decomposition}
\end{figure}

On peut définir l'espèce $\Arbo $, des \defn{arborescences}, comme solution de l'équation fonctionnelle
\begin{equation}\label{def_esparbo}
    \Arbo  =  \SingletonX\cdot (\Ens\circ \Arbo ).
 \end{equation}
Ceci correspond à décrire récursivement une arborescence comme étant constituée d'une racine (un point de $A$) à laquelle sont attachées des branches (qui sont chacune des arborescences sur des sous-ensembles) comme l'illustre la Figure~\ref{arbor}.
On a donc $ \Arbo (x)=x\,e^{\Arbo (x)}$.
L'espèce des arborescence intervient aussi dans la description de l'espèce des endofonctions. En effet, la décomposition de la Figure~\ref{endo_decomposition}  correspond à l'isomorphisme d'espèce  $\Endo =  \S(\Arbo )$. En passant aux séries, on trouve
ce qui donne
\begin{equation}\label{endo_ident}
   \Endo(x)=\frac{1}{ 1-\Arbo (x)}.
\end{equation}
Nous allons voir que cette identité permet de calculer le nombre d'arborescences à $n$ noeuds.
Introduisons d'abord deux autres opérations de base sur les espèces. 

\section{Dérivée et pointage.}\label{sec_derpoint}
\begin{figure}[ht]
\begin{center}
\scalebox{.75}{\includegraphics{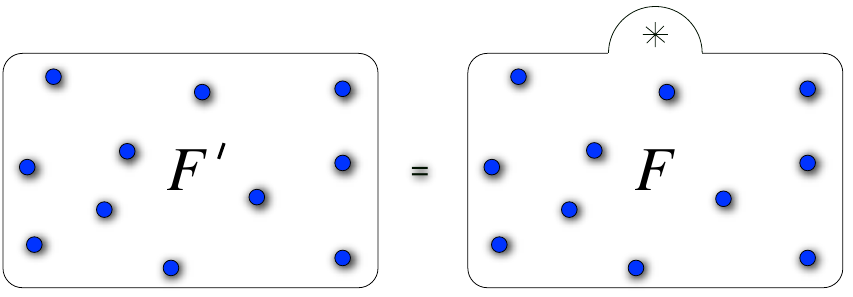}}
\end{center}
\vskip-10pt
\caption{Une structure typique d'espèce $\especeF'$.}
\label{deriv_fig}
\end{figure}

\begin{figure}[ht]\begin{center}
\scalebox{1}{\includegraphics{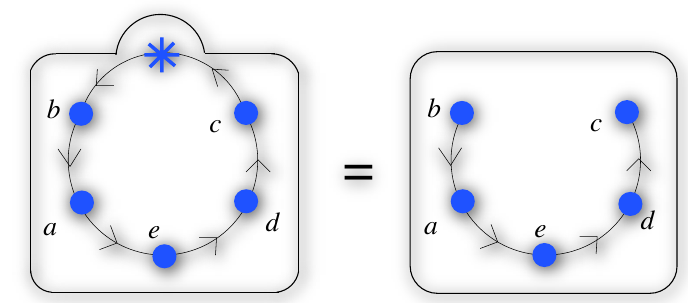}}
\end{center}
\vskip-10pt
\caption{La dérivée de l'espèce des cycles.}
\label{deriv_cyc}
\end{figure}

\begin{figure}[ht]\begin{center}
\scalebox{1}{\includegraphics{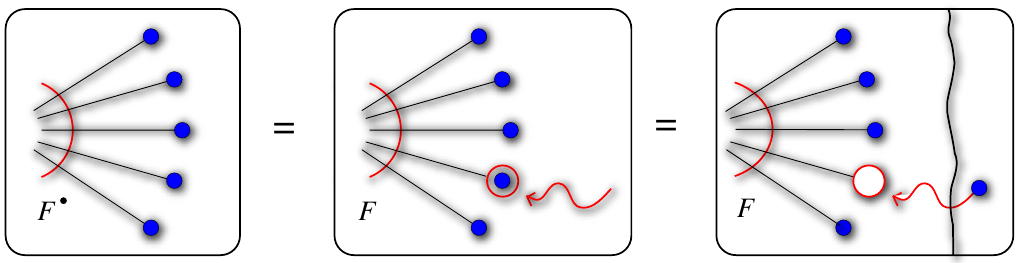}}
\end{center}
\vskip-10pt
\caption{Une structure typique d'espèce $\especeF\point$.}
\label{fig_pointage}
\end{figure}

\begin{figure}[ht]
\begin{center}
\scalebox{.8}{\includegraphics{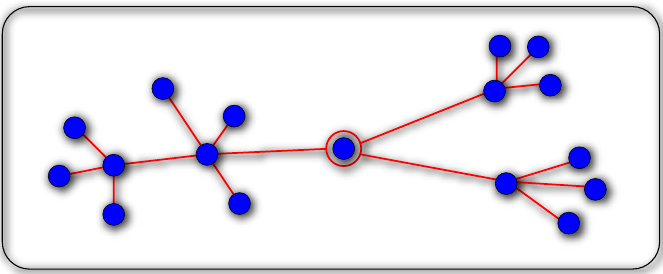}}
\end{center}
\vskip-10pt
\caption{Une arborescence \og est\fg\ un arbre pointé.}
\label{fig_arb_point}
\end{figure}

Les dernières opérations que l'on considère ici sont la \defn{dérivée} et le \defn{pointage} d'une espèce. Bien entendu, il y aura encore une fois compatibilité avec le passage aux séries.
Pour une espèce quelconque $\especeF$, les structures d'espèce $\especeF'$ sur un ensemble fini $A$ sont simplement les structures d'espèces $\especeF$ sur l'ensemble $A+\{*\}$. Une structure typique d'espèce $\especeF'$ peut donc se représenter à la Figure~\ref{deriv_fig}.
On illustre à la Figure~\ref{deriv_cyc} comment la dérivée de l'espèce des cycles peut s'identifier à l'espèce des ordre linéaires. 
En formule,
$\Cyc ' =  \Liste$, ce qui est exactement ce qu'on constate en passant aux séries génératrices correspondantes:
  $$\frac{d}{ dx}\, {\rm log}\frac{1}{ 1-x} = \frac{1}{ 1-x} .$$
Le \defn{pointage} d'une espèce est défini via la dérivée et la multiplication par $\SingletonX$, de la façon suivante. L'espèce $\especeF\point$ des $\especeF$ structures \defn{pointées} est $\especeF\point:=\SingletonX\cdot \especeF'$.
Autrement dit, une $\especeF\point$-structure sur $A$ est la donnée d'une $\especeF$-structure sur $A$ avec en plus le choix d'un point de $A$. Cette donnée s'identifie au choix de $a$ dans $A$ avec une $\especeF'$-structure sur $A\setminus \{a\}$ (voir la Figure~\ref{fig_pointage}). 
La justification de cette terminologie peut se voir par la Figure~\ref{fig_pointage}.
La Figure~\ref{fig_arb_point} montre qu'une arborescences est un {arbre}\footnote{Graphe simple connexe sans cycle.} pointé. 
On a donc 
\begin{equation}\label{arbre_arbo}
   \Arb\point=\SingletonX\cdot \Arbo ',
\end{equation}
 si $\Arb$ désigne l'\defn{espèce des arbres}. Comme en général, pour $\especeF(x)=\sum_{n=0}^\infty f_n\,{x^n/n!}$, on a
   $$\especeF\point(x)=\sum_{n=1}^\infty n\,f_n\, \frac{x^n}{ n!},$$
il découle de (\ref{arbre_arbo}) que le nombre d'arbres à $n$ noeuds est $a_n/n$ (avec $n\geq 1$), où $a_n$ désigne le nombre d'arborescences à $n$ noeuds. Nous allons maintenant calculer $a_n$.

\section{Vertébrés.}\label{sec_vertebre}
\begin{figure}[ht]
\begin{center}
\scalebox{1.2}{\includegraphics{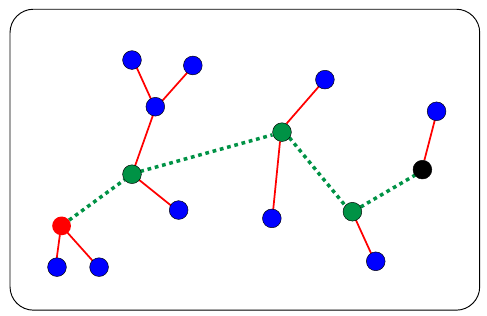}}
\begin{picture}(0,0)(0,0)
\put(-11,2.7){\tiny\rouge{Tête}}\put(-2.5,2.8){\tiny Queue}
\end{picture}
 \end{center}\vskip-15pt
\caption{Un vertébré.}
\label{fig_vertebre}
\end{figure}
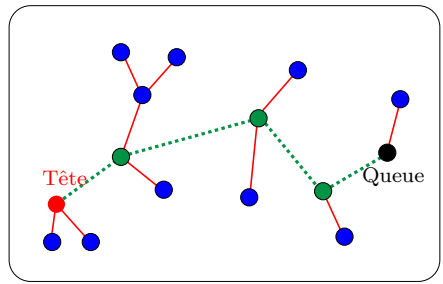

Notre énumération des arborescence passe par 
un argument élégant\footnote{Due à André Joyal.}, faisant appel à l'espèce des \og vertébrés\fg. Un \defn{vertébré} est une arborescence pointée, c'est-à-dire que l'espèce $\Vertebre$ des vertébrés est: $ \Vertebre  :=  \Arbo{\point}$.  On dit de la racine de l'arborescence que c'est la \defn{tête}  du vertébré. La \defn{queue}  est le noeud pointé dans l'arborescence. Il est tout-à-fait possible que ces deux sommets soient confondus.  L'unique chemin qui joint la tête à la queue est la  \defn{colonne vertébrale} du vertébré. Bien entendu, on dit des noeuds qui se trouvent sur cette colonne vertébrale, que ce sont les \defn{vertèbres} du vertébré. 
(voir Figure~\ref{fig_vertebre}).

\begin{figure}[ht]
\begin{center}
\scalebox{1}{\includegraphics{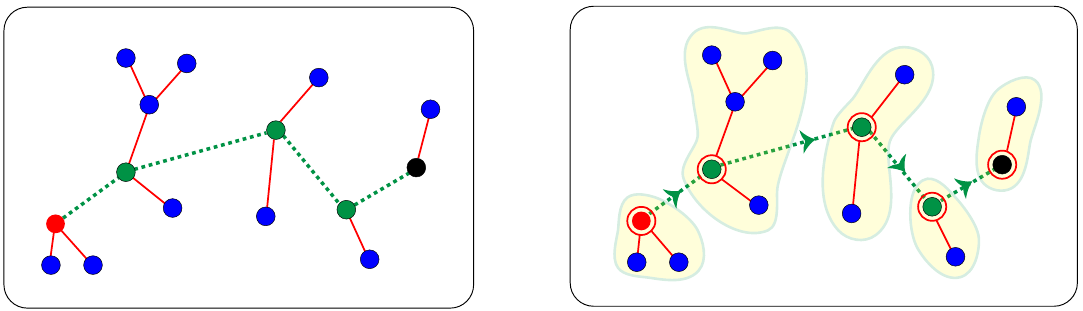}}
\begin{picture}(0,0)(0,0)
\put(0,0){\put(-22,2){\tiny\rouge{Tête}}\put(-14.5,2.3){\tiny Queue}}
\end{picture}
 \end{center}\vskip-15pt
\caption{Un vertébré transformé en liste d'arborescences.}
\label{fig_vertebre_decomp}
\end{figure}
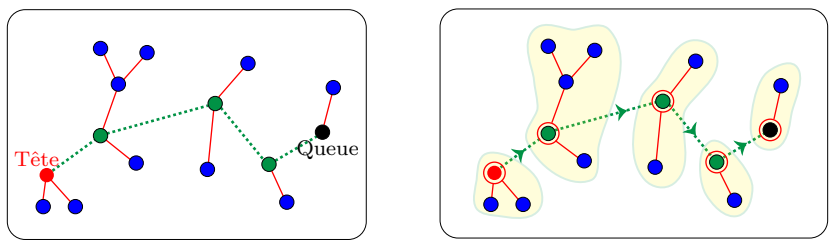

Dénotons $\nu_n$, le nombre de vertébrés sur un ensemble de cardinal  $n$.  Il est clair que
$ \nu_n  =  n\card{\Arbo[n]} $,
puisqu'il y a $n$ choix possibles pour  la queue.  
Calculons maintenant
$\nu_n$ d'une autre façon.  On oriente la colonne vertébrale de la tête vers la queue. Puisque à  
chaque vertèbre est attachée une arborescence, on obtient ainsi  
une liste non-vide d'arborescences disjointes, comme il est illustré à la
Figure~\ref{fig_vertebre_decomp}).

On met ainsi en évidence l'isomorphisme d'espèces
$\Vertebre   =  \Liste\nonvide (\Arbo)$, 
 de façon imagée. 
En passant aux séries génératrices associées, on trouve 
  $$\Vertebre(x)  = \frac{\Arbo(x)}{1-\Arbo(x)}.$$
Observons en passant que cette identité se déduit aussi directement de (\ref{def_esparbo}) par dérivation.
D'autre part, il découle directement de (\ref{endo_ident}) qu'on a aussi
 $$\Endo\nonvide(x) =  \frac{\Arbo(x)}{1-\Arbo(x)} =  \Vertebre(x).$$
On en déduit donc que  le nombre de vertébrés à $n$ sommets est égal au nombre d'endofonctions sur un ensemble à $n$ éléments, i.e.:  $\nu_n = n^n$. Finalement on a
\begin{prop}[Cayley 1889]
\bleu {Le nombre d'arborescences à $n$ noeud est $n^{n-1}$}.
\end{prop}

\section{Extensions et variantes}
Il est certain qu'il nous resterait à discuter bien d'autres aspects de la théorie des espèces. 
Nous n'avons qu'effleuré la théorie des espèces et décrit ses premières applications combinatoires. La théorie va beaucoup plus loin, et est toujours en développement. Pour en savoir plus, on peut consulter les articles originaux d'André Joyal~\cite{joyal, joyalanalytique}, ainsi que la monographie~\cite{BLL} qui présente l'essentiel des développements jusqu'en 1998. 
Cependant, pour donner une idée du spectre des notions combinatoires couvertes, on présente ci-dessous un rapide survol de certains de ces développements.  

Une partie un peu plus technique du développement de la théorie des espèces, qui met véritablement en lumière ses avantages par rapport aux autres approches proposées, permet d'intégrer la théorie de Pólya (voir~\cite{polya}) concernant l'énumération de structures à isomorphismes près (aussi considérées comme structures sur des objets non-étiquetés). Cette partie de la théorie des espèces se relie naturellement à la théorie de la représentation des groupes de permutations, et la théorie des fonctions symétriques (voir~\cite{bergeron}).

Plusieurs variantes de la théorie des espèces sont possibles, et même nécessaires. On a des espèces à plusieurs \og variables\fg\ (ou sortes), des espèces de structures avec des \og poids\fg, etc. Il faut aussi parfois adapter la théorie à certains contextes particuliers. 
Il est par exemple courant de travailler avec des ensembles qui viennent munis de façon explicite, ou implicite, d'un ordre sur les éléments. On peut vouloir considérer dans ces cas une théorie des espèces pour laquelle les structures s'élaborent sur des ensembles avec un ordre spécifié. Cette variante s'élabore facilement en imitant la théorie présentée dans ce texte, en tenant compte systématiquement de l'ordre donné sur les ensembles manipulés. Cela rend possible certaines constructions bien connues de la combinatoire classique, et facilite l'obtention des résultats d'énumération simples associés. On y perd cependant le lien avec la théorie de la représentation. 

En fait, la variante de la théorie des espèces avec laquelle on choisit de travailler dépend du type général des problèmes combinatoires considérés, tout comme l'espace de fonction dans lequel on travaille dépend du contexte considéré, qu'il soit algébrique, différentiel, analytique, etc.

La théorie des espèces trouve des applications non seulement en combinatoire et dans nombreux autres domaines des mathématiques (algèbre, théorie de la représentation, géométrie algébrique, etc.), mais aussi en physique \cite{faris,morton} et en informatique théorique \cite{carette}. La théorie a aussi été intégrée à certains environnements de calcul formels comme SAGE et Haskell.


\begin{thebibliography}{99}

\def\auteur#1{{\sc #1}}
\def\titreref#1{{\em #1}}
\def\vol#1{{\bf #1}}
\def\year#1{{\bf #1}}
\bibitem{bergeron} 
\auteur{F.~Bergeron},
\titreref{Algebraic Combinatorics and Coinvariant Spaces}, CMS Treatise in Mathematics, CMS and A.K.Peters,  2009.


\bibitem{BLL}
\auteur{F.~Bergeron, P.~Leroux, et G.~Labelle},
 \titreref{Combinatorial Species and Tree-Like Structures},
 Encyclopedia of Mathematics and its Applications, vol. \vol{67}, Cambridge University Press, \year{1998}.\\
 Une version remaniée du début de ce livre se retrouve ici:\\
  \url{http://bergeron.math.uqam.ca/Site/bergeron_anglais_files/livre_combinatoire_2.pdf}
  
 \bibitem{carette}
 \auteur{J.~Carette et G.~Uszkay},
 \titreref{Species: Making Analytic Functors Practical for Functional Programming}.
  
  
 \bibitem{faris}
 \auteur{W.~G.~Faris},
 \titreref{Combinatorics and Cluster Expansion},
Probab. Surveys \vol{7} (2010), 157--206.


 \bibitem{joyal}
 \auteur{A.~Joyal},
 \titreref{Une th\'eorie combinatoire des s\'eries formelles},
 Advances in Mathematics \vol{42} (1981), 1--82.
 
 \bibitem{joyalanalytique}	
 \auteur{A.~Joyal}, 
 \titreref{Foncteurs analytiques et esp\`eces de structures}, dans {\em Combinatoire \'Enum\'erative}, 
Springer Lecture Notes in Mathematics \vol{1234},
(1986),
  126--159. 
  
  \bibitem{jlabelle}
 \auteur{J. Labelle}, 
  \titreref{Applications diverses de la théorie combinatoire des espèces de structures}, 
  Annales des Sciences Math. du Québec, vol 7, 1983, no 1, p. 69--94. 

\bibitem{polya}	
\auteur{G.~P\'olya  and R.C.~Read}, 
\titreref{Combinatorial Enumeration of Groups, Graphs and Chemical Compounds},
Sprin\-ger--Verlag, Berlin, Heidelberg, and New York, 
 \year{1987}.

 \bibitem{morton}
 \auteur{Jeffrey Morton},
 \titreref{Categorified Algebra and Quantum Mechanics},
 Theory and Applications of Categories \vol{16} (2006), 785-854.


\end{thebibliography}
\end{document}